\documentclass[10pt,a4paper]{article}

\usepackage{amsmath,amssymb}
\usepackage{amsthm,latexsym}
\usepackage[labelformat=simple]{subfig}

\usepackage{bm,multirow}        
\usepackage{graphicx}  
\usepackage[software]{mymacros}

\graphicspath{{./fig/}}

\usepackage{hyperref}

\begin{document}

\title{Structure-Preserving Reduced Order Modeling of Non-Traditional Shallow Water Equation}

\author{S\"uleyman Y{\i}ld{\i}z \thanks{Institute of Applied Mathematics, Middle East Technical University, Ankara-Turkey
{yildiz.suleyman@metu.edu.tr}}, Murat Uzunca \thanks{Department of Mathematics, Sinop University, Turkey {muzunca@sinop.edu.tr}},
B\"ulent Karas\"ozen \thanks{Institute of Applied Mathematics \& Department of Mathematics, Middle East Technical University, Ankara-Turkey  {bulent@metu.edu.tr}}}

\date{}

\maketitle

\begin{abstract}

An energy preserving reduced order model is developed for the non-traditional shallow water equation (NTSWE) with full Coriolis force.
The NTSWE in the noncanonical Hamiltonian/Poisson form is discretized in space by finite differences. The resulting system of
ordinary differential equations is integrated in time by the energy preserving average vector field (AVF) method. The Poisson structure of the discretized NTSWE  exhibits a skew-symmetric matrix depending on the state variables.    An energy preserving, computationally efficient reduced order model
(ROM) is constructed by proper orthogonal decomposition with Galerkin projection. The nonlinearities are computed for the ROM efficiently by discrete empirical interpolation method. Preservation of the discrete energy and the discrete enstrophy are shown for the full order model,
and for the ROM which ensures the long term stability of the solutions.  The accuracy and computational
efficiency of the ROMs are shown by two numerical test problems.\\

\noindent\textbf{\textit{Keywords:}} Shallow water equation, model order reduction, Hamiltonian mechanics, finite difference methods, implicit time integrator.

\end{abstract}

\section{Introduction}

The shallow water equation (SWE) consists of a set of two-dimensional partial differential equations (PDEs) describing a thin inviscid fluid layer flowing over the topography in a rotating frame. SWE is a hyperbolic PDEs describing geophysical wave phenomena, e.g., the Kelvin and Rossby waves in the atmosphere and the oceans. SWEs are frequently used in large-scale geophysical flow prediction \cite{belanger2005data,cotter2012mixed}, investigation of baroclinic instability \cite{boss1996stability,vallis2017atmospheric}, and planetary flows \cite{warneford2014thermal}.  Energy and enstrophy are the most important conserved quantities of the SWEs, whereas the energy cascades to large scales whilst enstrophy cascades to small scales \cite{Cotter18,Stewart16}.

Real-time simulation of SWEs requires a large amount of computer memory and computing time. The reduced order models (ROMs) have emerged as a powerful approach to reduce the computational cost of evaluating large systems of PDEs like the SWE by constructing a low-dimensional linear reduced subspace, that approximately represents the solution to the system of PDEs with a significantly reduced computational cost. The solutions of the high fidelity full order model (FOM),  generated by space-time discretization of PDEs are projected usually on low dimensional reduced spaces using the proper orthogonal decomposition (POD), which is the widely used reduced order modeling technique. Applying POD Galerkin projection, the dominant POD modes of the PDEs are extracted from the snapshots of the FOM solutions. The computation of the FOM solutions is performed in the offline stage, whereas the reduced system from the low-dimensional subspace is solved in the online stage. The primary challenge in producing the low dimensional models of the high dimensional discretized PDEs is the efficient evaluation of the nonlinearities. The computational cost is reduced by sampling the nonlinear terms and interpolating, known as hyper-reduction techniques \cite{Astrid08,Barrault04,Carlberg13,chaturantabut10nmr,Nguyen08,Zimmermann16}.

The naive application of POD or DEIM may not preserve the geometric structures, like the symplecticness, energy preservation and passivity of Hamiltonian, Lagrangian and port Hamiltonian PDEs. The stability of reduced models over long-time integration and the structure-preserving properties has been recently investigated in the context of Lagrangian systems \cite{Carlberg15,Lall03}, and for port-Hamiltonian systems \cite{Chaturantabu16}. For linear and nonlinear Hamiltonian systems, the symplectic model reduction technique, proper symplectic decomposition (PSD) is constructed for Hamiltonian systems like linear wave equation, sine-Gordon equation, nonlinear Schr\"odinger equation to ensure long term stability of the reduced model \cite{Hesthaven16,Peng16}.  Recently the average vector field (AVF) method is used as a time integrator to construct reduced order models for Hamiltonian systems like  Korteweg-de Vries equation \cite{Gong17,Miyatake19} and nonlinear Schr\"odinger equation \cite{Karasozen18}.  Reduced order models for the SWEs are constructed in conservative form using POD-DEIM \cite{Lozovskiy17,Lozovskiy16}, in the $\beta$-plane by POD-DEIM and tensorial POD \cite{Stefanescu13,Stefanescu14}, by dynamic mode decomposition \cite{Bistrian15,Bistrian17}, the $f$-plane using POD \cite{Esfahanian09}. In these articles, the preservation of the energy and other conservative quantities in the reduced space are not discussed.

In this paper, we have constructed structure-preserving ROMs for the non-traditional shallow water equation (NTSWE) \cite{Dellar05,Stewart10,Stewart16} with the full  Coriolis force. Replacing the first order derivatives that appear in the NSTWE, a skew-gradient system, i.e. a non-canonical Hamiltonian system of ordinary differential equations (ODEs) is obtained. Time discretization of this system of non-canonical Hamiltonian system of ODEs by the AVF \cite{Cohen11} leads to FOM, which preserves the discrete Hamiltonian and Casimirs. The skew-symmetric structure of the full order skew-gradient system is preserved using the reduced order technique in \cite{Gong17,Karasozen18,Miyatake19}. The full order and reduced order NTSWE have state dependent skew-symmetric matrices, which does not allow separation of online and offline computation of the nonlinear terms.  Following \cite{Miyatake19} we have shown that the complexity of the ROM can be reduced for the POD and for the discrete empirical interpolation method (DEIM) \cite{chaturantabut10nmr}. The numerical results for two different representative examples  of the NTSWE confirm the structure preserving features like preserving the Hamiltonian (energy) and enstropy.  The efficiency of the ROMs are demonstrated by achieved speed-ups with the POD and DEIM over the FOM solutions.

The paper organized as follows. In Section~\ref{sec:swe}, the NTSWE is described in the Hamiltonian form. The structure preserving FOM in space and time is developed in Section~\ref{sec:fom}. The ROM with POD and DEIM are constructed in Section~\ref{sec:rom}. In Section~\ref{sec:num}, numerical results for two NTSWE examples are presented. The paper ends with some conclusions. 

\section{Shallow water equation}
\label{sec:swe}

Most of the models of the ocean and atmosphere include only the contribution to the Coriolis force from the component of the planetary rotation vector that is locally normal to geopotential surfaces when the vertical length scales are much smaller than the horizontal length scales.  This approach is known as  traditional approximation. However,  many atmospheric and oceanographic phenomena are substantially influenced by the non-traditional component of the Coriolis force \cite{stewart2013}, such as deep convection~\cite{marshall99}, Ekman spirals~\cite{leibovich1985}, and internal waves~\cite{gerkema2005}. The nondimensional NTSWE \cite{Dellar05,Stewart10,Stewart16} has the same structural form as the traditional SWE \cite{Salmon04} by distinguishing between the canonical velocities $\tilde{u}(x,y,t)$ and $\tilde{v}(x,y,t)$, and particle velocities $u(x,y,t)$ and $v(x,y,t)$
\begin{equation}\label{eq:ntswe}
\begin{aligned}
\frac{\partial \tilde{u}}{\partial t } & =  hqv - \frac{\partial \Phi }{\partial x},  \\
\frac{\partial \tilde{v}}{\partial t } & =  -hqu + \frac{\partial \Phi }{\partial y}, \\
\frac{\partial h}{\partial t } & =  - \frac{\partial }{\partial x} (hu) -  \frac{\partial }{\partial y} (hv)  , 
\end{aligned}
\end{equation}
where  $x$ and $y$ denote horizontal distances within a constant geopotential surface, and  $h(x,y,t)$ is the height field.
 The one-layer NTSWE \eqref{eq:ntswe}
describes an  inviscid fluid flowing over bottom topography at $z=h_b(x, y)$ in a frame rotating with angular velocity vector
 $ \bm{\Omega} =(\Omega^{(x)}, \Omega^{(y)}, \Omega^{(z)})$.  The orientation of the $x$ and $y$ axes are considered arbitrary with respect to North.
In traditional rotating and non-rotating SWEs, only the particle velocity components  appear. The canonical velocity components  are related to the canonical momentum per mass or to the depth average of particle velocities as
\begin{equation} \label{canpart}
\tilde{u}=u+2\Omega^{(y)}\left(h_b+\frac{1}{2}h\right), \quad \tilde{v}=v-2 \Omega^{(x)}\left(h_b+\frac{1}{2}h\right).
\end{equation}
The Bernoulli potential  $\Phi $ and potential vorticity $q$ are given by
\begin{align*}
\Phi &= \frac{1}{2} (u^2 + v^2) + g(h_b + h) + h(\Omega^{(x)} v  - \Omega^{(y)} u),   \\
 q &=\frac{1}{h}(2\Omega^{(z)}+\tilde{v}_x-\tilde{u}_y).
\end{align*}
 The non-traditional parameter is given as $ \delta=H/R_d $, where $ H $ represents the layer thickness scale and $ R_d $ is Rossby deformation radius, and $g$ denotes the gravitational acceleration \cite{Stewart10,Dellar05}.

The traditional SWE and NTSWE differ only by a function of the space alone, so their time derivatives are identical.
The non-rotating, traditional SWE \cite{Salmon88} and NTSWE  \eqref{eq:ntswe} have the same Hamiltonian structure and Poisson bracket \cite{Dellar05,Stewart10,Stewart16}
\begin{equation}\label{eq:ham_form}
\frac{\partial \tilde{z}}{\partial t} =
{ \mathcal J}(\tilde{z})\dfrac{\delta \mathcal{H}}{\delta z} 
= \begin{pmatrix}
 0& q&-\partial_x \\
 -q& 0&-\partial_y \\
 -\partial_x&-\partial_y  & 0
 \end{pmatrix}
 \begin{pmatrix}
hu \\
hv\\
\Phi
 \end{pmatrix},
\end{equation}
where $z=(u,v,h)$ and $\tilde{z} = (\tilde{u},\tilde{v},h)$.
The Hamiltonian or the energy of \eqref{eq:ntswe} is given in terms of particle velocity components by
\begin{equation} \label{eq:hamilton1}
 {\mathcal H}(z) =  \iint  \bigg\{\frac{1}{2}h(u^2 + v^2) +gh\left (h_b + \frac{1}{2}h\right)
\bigg\}d{\mathbf x},
\end{equation}
over a periodic domain.
We remark that the Hamiltonian \eqref{eq:hamilton1} is treated as a function of the canonical velocity components $\tilde{u}$ and $\tilde{v}$ and the layer thickness using the relations \eqref{canpart}.

The non-canonical Hamiltonian form of NTSWE   \eqref{eq:ham_form} is determined by the skew-adjoint Poisson bracket of two functionals $\mathcal{A}$ and $\mathcal{B}$ \cite{Lynch02,Salmon04} as
\begin{equation} \label{bracket}
\{ \mathcal{A},\mathcal{B}  \} = \iint \left(q\frac{\delta(\mathcal{A},\mathcal{B})}{\delta (\tilde{u},\tilde{v}) } - \frac{\delta \mathcal{A} }{\delta \tilde{\bm{\upsilon}}}\cdot \nabla \frac{\delta \mathcal{B}}{\delta h} + \frac{\delta \mathcal{B}}{\delta \tilde{\bm{\upsilon}}}\cdot \nabla \frac{\delta \mathcal{A}}{\delta h} \right) d {\mathbf x},
\end{equation}
where $\tilde{\bm{\upsilon}} = (\tilde{u},\tilde{v})$.
The functional Jacobian is given by
\begin{equation*}
\frac{\delta(\mathcal{A},\mathcal{B})}{\delta(\tilde{u},\tilde{v})} = \frac{\delta \mathcal{A}}{\delta \tilde{u}}\frac{\delta
 \mathcal{B}}{\delta \tilde{v} } - \frac{\delta \mathcal{B}}{\delta \tilde{u}}\frac{\delta \mathcal{A}}{\delta \tilde{v}}.
\end{equation*}
The Poisson bracket  \eqref{bracket} is related to the skew-symmetric Poisson matrix ${\mathcal J}$ as $\{\mathcal{A},\mathcal{B}\}=\{\mathcal{A},{\mathcal J}\mathcal{B}\}$.   Although
the matrix  ${\mathcal J}$ in \eqref{eq:ham_form} is not skew-symmetric, the skew-symmetry of the Poisson
bracket appears after integrations by parts \cite{Lynch02}, and the Poisson bracket satisfies the Jacobi identity
\begin{equation*}
\{ \mathcal{A},\{\mathcal{B},\mathcal{D}\}\} + \{ \mathcal{B},\{\mathcal{D},\mathcal{A}\}\} + \{ \mathcal{A},\{\mathcal{B},\mathcal{D}\}\} = 0,
\end{equation*}
for any three functionals $\mathcal{A}$, $\mathcal{B}$ and $\mathcal{D}$.
The preservation of the Hamiltonian follows from the antisymmetry of the Poisson bracket \eqref{bracket}
\begin{equation*}
\frac{d{\mathcal H}}{dt}= \{{\mathcal H}, {\mathcal H} \} = 0.
\end{equation*}
Besides the Hamiltonian, there are other conserved quantities in form of Casimirs
$$
\mathcal{C} = \iint  hG(q)d{\mathbf x},
$$
where $G$ is an arbitrary function of the potential vorticity
$q$.
The Casimirs are  additional constants of motion which commute with any functional $\mathcal{A}$, i.e  the Poisson bracket vanishes.
Important special cases are the potential enstrophy
$$
{\mathcal Z } =  \frac{1}{2}  \int \int hq^2 d{\mathbf x} =  \frac{1}{2}  \int \int \frac{1}{h}\left (\Omega^{(z)} +\frac{\partial \tilde{v}}{\partial x}   - \frac{\partial \tilde{u}}{\partial y}   \right )^2 d{\mathbf x},
$$
the mass ${\mathcal M } =\iint hd{\mathbf x}$, and the vorticity ${\mathcal V } =\iint h q d{\mathbf x}$.

\section{Full order model}
\label{sec:fom}

The NTSWE \eqref{eq:ntswe} is discretized by finite differences on a uniform grid in the spatial domain $(a,b)\times(c,d)$ with the  nodes ${\mathbf x}_{ij} = (x_i,y_j )^T$, where $x_i=a+(i-1)\Delta x$ and $y_j=c+(j-1)\Delta y$, $i=1,\ldots, N_x+1$, $j=1,\ldots, N_y+1$, and then
discretized in space canonical and particle velocity components and height are given by
\begin{equation}\label{solvec}
\begin{aligned}
{\mathbf u}(t) &= (u_{11}(t),\ldots , u_{1N_y}(t),u_{21}(t),\ldots , u_{2N_y}(t), \ldots , u_{N_xN_y}(t))^T,\\
{\mathbf v}(t) &= (v_{11}(t),\ldots , v_{1N_y}(t),v_{21}(t),\ldots , v_{2N_y}(t), \ldots , v_{N_xN_y}(t))^T,\\
\tilde{\mathbf u}(t) &= (\tilde{u}_{11}(t),\ldots , \tilde{u}_{1N_y}(t),\tilde{u}_{21}(t),\ldots , \tilde{u}_{2N_y}(t), \ldots , \tilde{u}_{N_xN_y}(t))^T,\\
\tilde{\mathbf v}(t) &= (\tilde{v}_{11}(t),\ldots , \tilde{v}_{1N_y}(t),\tilde{v}_{21}(t),\ldots , \tilde{v}_{2N_y}(t), \ldots , \tilde{v}_{N_xN_y}(t))^T,\\
{\mathbf h}(t) &= (h_{11}(t),\ldots , h_{1N_y}(t),h_{21}(t),\ldots , h_{2N_y}(t), \ldots , h_{N_xN_y}(t))^T.
\end{aligned}
\end{equation}
where for $w=u,v,\tilde{u},\tilde{v},h$, $w_{ij}(t)$ denotes the approximation of $w({\mathbf x},t)$ at the grid nodes ${\mathbf x}_{ij}$ at time $t$, $i=1,\ldots, N_x$, $j=1,\ldots, N_y$.
We note that the degree of freedom is given by $N=N_xN_y$ because of the periodic boundary conditions, i.e., the most right and the most top grid nodes are not included. 
Throughout the paper, we do not explicitly represent the time dependency of the semi-discrete solutions for simplicity, and we write ${\mathbf u}$, ${\mathbf v},\tilde{\mathbf u}$, $\tilde{\mathbf v}$ and ${\mathbf h}$.
The semi-discrete vector for the solution vectors are  defined by ${\mathbf z}=({\mathbf u},{\mathbf v},{\mathbf h})\in\mathbb{R}^{3N}$ and
$\tilde{\mathbf z}=(\tilde{\mathbf u},\tilde{\mathbf v},{\mathbf h})\in\mathbb{R}^{3N}$.

For the approximation of the first order partial derivative terms, we use one dimensional central finite differences to the first order derivative terms in either $x$ and $y$ direction, and we extend them to two dimensions utilizing the Kronecker product.
For a positive integer $s$, let $\widetilde{D}_s$ denotes the matrix related to the  central finite differences  to the first order ordinary differential operator under periodic boundary conditions
$$
 \widetilde{D}_s=
\begin{pmatrix}
 0& 1&  & &-1 \\
-1& 0&1 & &   \\
  & \ddots  & \ddots  &\ddots  &   \\
  &  &-1&0 &1 \\
 1&  &  & 1&0
\end{pmatrix} \in \mathbb{R}^{s\times s}.
$$
Then, on the two dimensional mesh, the central finite difference matrices corresponding to the first order partial derivative operators $\partial_x$ and $\partial_y$ are given respectively by
$$
D_x=\frac{1}{2\Delta x}\widetilde{D}_{N_x}\otimes I_{N_y}\in\mathbb{R}^{N\times N} \; ,  \quad D_y=\frac{1}{2\Delta y}I_{N_x}\otimes \widetilde{D}_{N_y}\in\mathbb{R}^{N\times N},
$$
where $\otimes$ denotes the Kronecker product, and
 $I_{N_x}$ and $I_{N_y}$ are the identity matrices of size $N_x$ and $N_y$, respectively.

We further partition the time interval $[0,T]$ into $N_t$ uniform intervals with the step size $\Delta t=T/N_t$ as $0=t_0<t_1<\ldots <t_{N_t}=T$, and $t_k=k\Delta t$, $k=0,1,\ldots ,N_t$. Then, we denote by $\tilde{\mathbf u}^k=\tilde{\mathbf u}(t_k)$, $\tilde{\mathbf v}^k=\tilde{\mathbf v}(t_k)$ and ${\mathbf h}^k={\mathbf h}(t_k)$ the full discrete solution vectors at time $t_k$. Similar setting is used for the other components, as well.

The full discrete form of the energy and the enstrophy at a time instance $t_k$ are given as
\begin{align} 
H^k &= \sum_{i=1}^N \bigg\{ \frac{1}{2}{\mathbf h}^k_i\left(({\mathbf u}^k_i)^2  + ({\mathbf v}^k_i)^2\right) + g{\mathbf h}^k_i\left(({\mathbf h}_b)_i + \frac{1}{2}{\mathbf h}^k_i\right)\bigg\}\Delta x\Delta y, \label{dener} \\
Z^k &= \frac{1}{2} \sum_{i=1}^N \frac{\left((D_x{\tilde{\mathbf v}}^k)_i - (D_y{\tilde{\mathbf u}}^k)_i + \Omega^{(z)}\right)^2}{{\mathbf h}^k_i}\Delta x\Delta y. \nonumber
\end{align}

The semi-discrete formulation of the NTSWE \eqref{eq:ntswe} leads to a $3N$ dimensional system of Hamiltonian ODEs in skew-gradient form
\begin{align}\label{semi-dis-poiss}
\dfrac{d \tilde{\mathbf z}}{d t}&= J(\tilde{\mathbf z})\nabla_{{\mathbf z}} H({\mathbf z})
=
\begin{pmatrix}
0& {\mathbf q}^d&-D_x \\
-{\mathbf q}^d& 0&-D_y \\
-D_x&-D_y  & 0
\end{pmatrix}
\begin{pmatrix}
{\mathbf u} \circ {\mathbf h} \\
{\mathbf v}\circ {\mathbf h} \\
\bm{\Phi}
\end{pmatrix}, 
\end{align}
with the discrete Bernoulli potential
$$
\bm{\Phi} = \frac{1}{2} ({\mathbf u}\circ {\mathbf u} +{\mathbf v}\circ {\mathbf v})+g({\mathbf h}+\bh_b) +{\mathbf h}
\left(\Omega^{(x)}{\mathbf v} -\Omega^{(y)}{\mathbf u}\right),
$$
where $ \circ $ denotes element-wise or Hadamard product. The matrix  ${\mathbf q}^d\in \mathbb{R}^{N\times N}$ is the diagonal matrix with the diagonal elements ${\mathbf q}^d_{ii}={\mathbf q}_i$ where ${\mathbf q}$ is the semi-discrete vector of the potential vorticity $q$, $i=1,\ldots ,N$.

For time integration we use the Poisson structure preserving average vector field method (AVF). The AVF method  preserves higher order polynomial Hamiltonians \cite{Cohen11}, including the cubic Hamiltonian ${\mathcal H}$ of the NTSWE \eqref{eq:ntswe}.
Quadratic Casimirs function like mass and circulation are preserved exactly by AVF method. But higher-order polynomial Casimirs like the enstrophy (cubic) can  not  be preserved.  Practical implementation of the AVF method requires the evaluation of the integral on the right-hand side \eqref{swe_avf}. Since the Hamiltonian ${\mathcal H}$ and the discrete form of the Casimirs, potential enstrophy, mass and circulation are polynomial, they can be exactly integrated with a Gaussian quadrature rule of the appropriate degree.
The AVF method is used with finite element discretization of the rotational SWE  \cite{Cotter18,Cotter19} and for thermal SWE \cite{Eldred19} in Poisson form. After time integration of the semidiscrete NTSWE  \eqref{semi-dis-poiss} by the AVF integrator, the full discrete problem reads as: for $k=0,1,\ldots , N_{t}-1$, given $\tilde{\mathbf z}^{k}$ find $\tilde{\mathbf z}^{k+1}$ satisfying
\begin{align}\label{swe_avf}
 \tilde{\mathbf z}^{k+1}&=\tilde{\mathbf z}^{k}+ \Delta t J\left (\dfrac{\tilde{\mathbf z}^{k+1}+\tilde{\mathbf z}^{k}}{2}\right )\int_{0}^{1} \nabla_{\mathbf z} H(\xi({\mathbf z}^{k+1}-{\mathbf z}^{k})+{\mathbf z}^{k})d \xi.
\end{align}

\section{Reduced order model}
\label{sec:rom}

In this section, we construct ROMs that preserve the skew-gradient structure of the NTSWE  \eqref{semi-dis-poiss} and consequently the discrete Hamiltonian \eqref{dener}.  Because the NTSWE is a non-canonical Hamiltonian PDE with a state dependent Poisson structure, a straightforward application of the POD will not preserve the skew-gradient structure of the NTSWE  \eqref{semi-dis-poiss} in reduced form. Energy  preserving POD reduced systems are constructed for Hamiltonian systems with constant skew-symmetric matrices like the Korteweg de Vries equation  \cite{Gong17,Miyatake19} and nonlinear Schr\"odinger equation \cite{Karasozen18}. The approach in \cite{Gong17} can  be applied to skew-gradient systems with state dependent skew-symmetric structure as the NTSWE \eqref{semi-dis-poiss}. 
We show that the sate dependent skew-symmetric matrix in \eqref{semi-dis-poiss} can be evaluated efficiently in the online stage independent of the full dimension $N$.

The POD basis are computed through the mean subtracted snapshot matrices $S_{\tilde {u}}$, $S_{\tilde {v}}$ and $S_h$, constructed by the solutions of the full discrete high fidelity model \eqref{swe_avf}
\begin{align*}
S_{\tilde {u}} &= \left(\tilde{\mathbf u}^1 - \overline{\tilde{\mathbf u}},  \cdots, \tilde{\mathbf u}^{N_t} - \overline{\tilde{\mathbf u}} \right) \in\mathbb{R}^{N\times N_t}, \\
S_{\tilde {v}} &= \left(\tilde{\mathbf v}^1 - \overline{\tilde{\mathbf v}}, \cdots, \tilde{\mathbf v}^{N_t} - \overline{\tilde{\mathbf v}} \right)   \in\mathbb{R}^{N\times N_t},  \\
S_h &= \left({\mathbf h}^1 - \overline{\mathbf h}, \cdots, {\mathbf h}^{N_t}  - \overline{\mathbf h}\right) \in\mathbb{R}^{N\times N_t},
\end{align*}
where $\overline{\tilde{\mathbf u}}$, $\overline{\tilde{\mathbf v}}$, $\overline{\mathbf{h}}\in\mathbb{R}^{N}$  denote the time averaged mean of the solutions
$$
\overline{\tilde{\mathbf u}} = \frac{1}{N_t}\sum_{k=0}^{N_t} \tilde{\mathbf u}^k\; , \quad \overline{\tilde{\mathbf v}} = \frac{1}{N_t}\sum_{k=0}^{N_t} \tilde{\mathbf v}^k\; , \quad \overline{\mathbf h} = \frac{1}{N_t}\sum_{k=0}^{N_t} { \mathbf h}^k.
$$
The mean-subtracted ROMs is used frequent in fluid dynamics, and it guarantees that ROM solution would
satisfy the same boundary conditions as the FOM.

The POD modes are computed by applying singular value decomposition (SVD) to the snapshot matrices
\begin{equation*}
S_{\tilde{u}}= W_{\tilde{u}} \Sigma_{\tilde{u}} U_{\tilde{u}}^T\; , \quad S_{\tilde{v}}= W_{\tilde{v}} \Sigma_{\tilde{v}} U_{\tilde{v}}^T \; , \quad S_h= W_h \Sigma_h U_h^T ,
\end{equation*}
where for $i=\tilde{u},\tilde{v},h$, the columns of the orthonormal matrices $W_i \in \mathbb{R}^{ N\times N_t}$ and
$U_i\in \mathbb{R}^{ N_t\times N_t}$ are the left and right singular vectors of the snapshot matrices $S_i$, respectively, and the diagonal matrix $\Sigma_i\in \mathbb{R}^{ N_t\times N_t}$ contains the singular values $\sigma_{i,1} \ge \sigma_{i,2} \ge \cdots \ge \sigma_{i,N_t}\geq 0$.
Then, the matrix $V_{i,n}\in\mathbb{R}^{N\times n}$ of rank $n$ POD modes consists of the first $n$ left singular vectors from $W_i$ corresponding to the $n$ largest singular values, which satisfies
the following least squares error
$$
\min_{V_{i,n}\in \mathbb{R}^{ N\times n}}||S_i-V_{i,n}V_{i,n}^TS_i ||_F^2 = \sum_{j=n+1}^{N_t} \sigma_{i,j}^2\; , \quad i=\tilde{u},\tilde{v},h,
$$
where $\|\cdot\|_F$ is the Frobenius norm. Moreover, we have the reduced approximations
\begin{equation}\label{relz}
\tilde{\mathbf u} \approx  \overline{\tilde{\mathbf u}} + V_{\tilde{\mathbf u},n}\tilde{\mathbf u}_r, \quad \tilde{\mathbf v} \approx  \overline{\tilde{\mathbf v}} + V_{\tilde{\mathbf v},n}\tilde{\mathbf u}_r,
\quad {\mathbf  h} \approx \overline{\mathbf h} + V_{h,n}  {\mathbf h}_r,
\end{equation}
where the reduced (coefficient) vectors $\tilde{\mathbf u}_r$, $\tilde{\mathbf u}_r$ and ${\mathbf h}_r$ are the solutions of the following ROM of \eqref{semi-dis-poiss}
\begin{align} \label{galpod1}
\frac{d}{dt}\tilde{\mathbf z}_r = V_{z,n}^T J(\tilde{\mathbf z})\nabla_{\mathbf z} H({\mathbf z}),
\end{align}
where $\tilde{\mathbf z}_r=(\tilde{\mathbf u}_r,\tilde{\mathbf v}_r,{\mathbf h}_r)$, and the components of the vector $\tilde{\mathbf z}=(\tilde{\mathbf u},\tilde{\mathbf v},{\mathbf h})$ are given as in \eqref{relz}. The block diagonal matrix $V_{z,n}$ contains the matrix of POD modes for each solution component given by
\begin{equation*}
V_{z,n}=
\begin{pmatrix}
V_{\tilde{\mathbf u},n} & &\\
& V_{\tilde{\mathbf v},n}&\\
& & V_{{\mathbf h},n}
\end{pmatrix}\in\mathbb{R}^{3N\times 3n}.
\end{equation*}

The ROM \eqref{galpod1} is not a skew-gradient system.
A reduced skew-gradient system is obtained formally by inserting $V_{z,n}V_{z,n}^T$ between $J(\tilde{\mathbf z})$ and $\nabla_{\mathbf z} H({\mathbf z})$  \cite{Gong17}, leading to the ROM
\begin{align} \label{galpod}
\frac{d}{dt}\tilde{\mathbf z}_r = J_r(\tilde{\mathbf z})\nabla_{\mathbf{z}_r} H({\mathbf z}),
\end{align}
where $ J_r(\tilde{\mathbf z})= V_{z,n}^T J( \tilde{\mathbf z})V_{z,n}$ and $\nabla_{\mathbf{z}_r} H({\mathbf z})= V_{z,n}^T\nabla_{\mathbf{z}} H({\mathbf z})$.  The reduced order NTSWE  \eqref{galpod} is also solve by the AVF.

The reduced NTSWE \eqref{galpod} can be written explicitly as
\begin{align}\label{skewred}
\frac{d}{dt}\tilde{\mathbf z}_r &=\begin{pmatrix}
0&V_{u,n}^T\bq^{d}V_{v,n}&-V_{u,n}^TD_xV_{h,n}\\
-V_{v,n}^T\bq^{d}V_{u,n}&0 & -V_{v,n}^TD_yV_{h,n}\\
-V_{h,n}^TD_xV_{u,n}     & -V_{h,n}^TD_yV_{v,n}&0
\end{pmatrix}
V_{z,n}^T\nabla_{\mathbf{z}} H(\mathbf{z}).
\end{align}
The reduced system \eqref{skewred} has constant matrices which can be  precomputed in offline stage  whereas the matrices $ V_{u,n}^T\bq^{d}V_{v,n} $ and $ V_{v,n}^T\bq^{d}V_{u,n} $ should be computed in online stage depending on the full order system. Exploiting the diagonal structure of $ \bq^{d} $ the computational complexity of evaluating the state dependent skew-symmetric matrix in \eqref{skewred} can be reduced  similar to the skew-gradient systems with constant skew-symmetric matrices  as in \cite{Miyatake19}.
Let $ \text{vec}(\cdot) $ denotes vectorization of a matrix.
For any $  A \in \mathbb{R}^{m\times n}  $ and $ B \in \mathbb{R}^{n\times p}$
\begin{align*}
\text{vec}(AB) =(I_p\otimes A)\text{vec(B)} =(B^\top\otimes I_m)\text{vec(A)}.
\end{align*}
Thus, for a diagonal matrix $  D \in \mathbb{R}^{n\times n}  $ and $ V \in \mathbb{R}^{n\times r} $
\begin{align*}
\text{vec}(V^\top DV)&=(I_r \otimes V^\top)\text{vec}(DV)\\
&=(I_r \otimes V^\top)(V^\top \otimes I_n)\text{vec}(D)\\
&=(V\otimes V)^\top\text{vec}(D)\\
&=(V\otimes V)^\top\tilde{M}^\top\tilde{D}\\
&=\begin{pmatrix}
V(1,:)\otimes V(1,:)\\
\vdots\\
V(n,:)\otimes V(n,:)
\end{pmatrix}^\top
\tilde{D},
\end{align*}
where $M$ is an operator satisfying $ M(\bz\otimes\bz)=\bz\circ\bz $ and $ \tilde{D}= [D_{11},D_{22},\ldots,D_{nn}]^T $. Using the above result, the computational complexity of the matrix products $ V_{u,n}^T\bq^{d}V_{v,n} $ and $ V_{v,n}^T\bq^{d}V_{u,n}$ is reduced from $ \mathcal{O}(n\cdot N(n+N)) $ to $  \mathcal{O}(n^2\cdot N) $.

Due to nonlinear terms, the computation of the reduced system still scales with the dimension $N$ of the FOM. This can be reduced  by applying the hyper-reduction technique such as  DEIM \cite{chaturantabut10nmr}.
The ROM \eqref{galpod} can be rewritten as a nonlinear ODE system  of the form
\begin{equation}\label{rom-semi}
\frac{d}{dt}\tilde{\mathbf z}_r = V_{z,n}^T F(\tilde{\mathbf z})=\begin{pmatrix}
V_{u,n}^TF_1(\tilde{\mathbf z}) \\
V_{v,n}^TF_2(\tilde{\mathbf z}) \\
V_{h,n}^TF_3(\tilde{\mathbf z})
\end{pmatrix}.
\end{equation}
The DEIM is applied by
sampling the nonlinearity $F(\cdot)$ and then interpolating with hyper-reduction. To obtain the DEIM basis, we form the snapshot matrices defined by
$$
G_i = ( F_i^1,F_i^2,\cdots , F_i^{N_t})\in\mathbb{R}^{N\times N_t}, \quad i=1,2,3,
$$
where $F_i^k=F_i(\tilde{\mathbf z}^k)$ denotes the $i$-th component of the nonlinearity $F(\tilde{\mathbf z})$ in \eqref{rom-semi} at time $t_k$ computed by using the FOM solution vector $\tilde{\mathbf z}$, $k=1,\ldots , N_t$. Then, we can approximate each $F_i(\tilde{\mathbf z})$ in the column space of the snapshot matrices $G_i$. We first apply POD to the snapshot matrices $G_i$ and find the basis matrices $V_{F_i,m}\in\mathbb{R}^{N\times m}$ whose columns are the basis vectors spanning the column space of the snapshot matrices $G_i$.  We apply the DEIM algorithm \cite{chaturantabut10nmr} to find a projection matrix $P_i\in\mathbb{R}^{N\times m}$
\begin{equation*}
F_i(\tilde{\mathbf z}) \approx V_{F_i,m}(P_i^TV_{F_i,m})^{-1} P_i^TF_i(\tilde{\mathbf z}),
\end{equation*}
and then we get the DEIM approximation to the reduced nonlinearities in \eqref{rom-semi} as
\begin{equation*}
V_{u,n}^TF_1(\tilde{\mathbf z}) \approx \mathcal{V}_{u,1} (P_1^TF_1(\tilde{\mathbf z})), \quad V_{v,n}^TF_2(\tilde{\mathbf z}) \approx \mathcal{V}_{v,2} (P_2^TF_2(\tilde{\mathbf z})), \quad V_{h,n}^TF_3(\tilde{\mathbf z}) \approx \mathcal{V}_{h,3} (P_3^TF_3(\tilde{\mathbf z})),
\end{equation*}
where
\begin{equation*}
\mathcal{V}_{u,1} = V_{u,n}^TV_{F_1,m}(P_1^TV_{F_1,m})^{-1}, \; \mathcal{V}_{v,2} = V_{v,n}^TV_{F_2,m}(P_2^TV_{F_2,m})^{-1}, \;  \mathcal{V}_{h,3} = V_{h,n}^TV_{F_3,m}(P_3^TV_{F_3,m})^{-1}
\end{equation*}
 are all the matrices of size $n\times m$, and they are precomputed in the offline stage.
Using the DEIM approximations, the ROM \eqref{rom-semi} becomes

\begin{equation*}
\frac{d}{dt}\tilde{\mathbf z}_r = \begin{pmatrix}
\mathcal{V}_{u,1}F_{r,1}(\tilde{\mathbf z}) \\
\mathcal{V}_{v,2}F_{r,2}(\tilde{\mathbf z}) \\
\mathcal{V}_{h,3}F_{r,3}(\tilde{\mathbf z})
\end{pmatrix},
\end{equation*}
where the reduced nonlinearities $F_{r,i}(\tilde{\mathbf z})=P_i^TF_i(\tilde{\mathbf z})$ are computed by considering just $m\ll N$ entries of the nonlinearities $F_i(\tilde{\mathbf z})$ among $N$ entries, $i=1,2,3$.

\section{Numerical results}
\label{sec:num}

In this section we present two numerical examples to demonstrate the efficiency of the ROMs. We consider the propagation of the inertia-gravity waves by Coriolis force, known as  geostrophic adjustment \cite{Stewart16}, and the shear instability in the form of roll-up of an unstable shear layer, known as barotropic instability  \cite{Stewart16}.
For numerical simulations, we consider the nondimensional form of the NTSWE \eqref{eq:ntswe}  with the setting
$$
x = R_d\hat{x},\quad y = R_d\hat{y},\quad u=c\hat{u},\quad v= c\hat{v}, \quad h =H\hat{h}, \quad h_b= H\hat{h}_b,
$$
$$
\left(\Omega^{(x)},\Omega^{(y)},\Omega^{(z)}\right)
= \Omega \left(\hat{\Omega}^{(x)},\hat{\Omega}^{(y)},\hat{\Omega}^{(z)}\right),
$$
where a $\hat{ }$ denotes a dimensionless variable, and $\Omega$ is planetary rotation rate to construct the gravity wave speed $c$
$$
c = \sqrt{gH}, \quad R_d = \frac{c}{2\Omega}, \quad \delta = \frac{H}{R_d} = \frac{2\Omega H}{c}.
$$
The parameters are taken following \cite{Stewart16} as $H=1000$ m, $\Omega\approx 7.3 \times 10^{-5}$ rad s $^{-1}$, $g=10^{-3}$ms$^{-2}$. The dimensionless components of the rotation vector at latitude $\phi$ are taken as
$$
\hat{\Omega}^{(x)} =0, \quad \hat{\Omega}^{(y)} =\cos(\phi), \quad \hat{\Omega}^{(z)} =\sin(\phi),
$$
where  we set $\phi=\pi/4$ in the numerical experiments. In all examples,
the spatial and temporal mesh sizes are taken as $\Delta x = 0.1$ and $\Delta t=  0.1$, respectively.

In order to determine the numbers $n$ and $m$ of the POD and DEIM modes, respectively, we use the so-called relative cumulative energy criteria for a desired number $p=m,n$
\begin{equation} \label{energy_criteria}
\min_{p}\frac{\sum_{j=1}^p \sigma_{j}^2}{\sum_{j=1}^{N_t} \sigma_{j}^2  } > 1 - \kappa,
\end{equation}
where $\kappa $ is a user-specified tolerance. In our simulations, we set $\kappa = 10^{-3}$ and $\kappa = 10^{-5}$
to catch at least $99.9 \%$ and $99.999 \%$ of data information for POD and DEIM modes, respectively. We take  the same number of modes for each state variable.

The error between a discrete FOM solution and a discrete reduced approximation (FOM-ROM error) are measured for the components ${\mathbf w}= \tilde{\mathbf u}, \tilde{\mathbf v},{\mathbf h}$  using the following time averaged relative $L_2$ errors
\begin{align*}
\|\mathbf{w}-\widehat{\mathbf w}\|_{rel}=\frac{1}{N_t}\sum_{k=1}^{N_t}\frac{\|{\mathbf w}^k-\widehat{\mathbf w}^k\|_{L^2}}{\|{\mathbf w}^k\|_{L^2}}, \quad  \|{\mathbf w}^k\|_{L^2}^2=\sum_{i=1}^N({\mathbf w}^k_i)^2\Delta x\Delta y,
\end{align*}
where $\widehat{\mathbf w}=\overline{\mathbf w} + V_{{\mathbf w},n}{\mathbf w}_r$ denotes the reduced approximation to ${\mathbf w}$.
All simulations are performed on a machine with Intel$^{\circledR}$
Core$^{{\mathrm TM}}$ i7 2.5 GHz 64 bit CPU, 8 GB RAM, Windows 10, using 64 bit MatLab R2014.

\subsection{Single-layer geostrophic adjustment}
\label{ex1}

We consider the NTSWE on the periodic spatial domain $[-5, 5]^2$ and on the time interval $[0,100]$ \cite{Stewart16}.
The initial conditions are prescribed in form of a motionless layer with an upward bulge of the height field
\begin{align*}
&h(x, y, 0) = 1+\frac{1}{2} \exp\left[-\left(\frac{4x}{5}\right)^2-\left(\frac{4y}{5}\right)^2\right], \nonumber\\
&u(x, y, 0) = 0,\\
&v(x, y, 0) = 0.\nonumber
\end{align*}
The inertia-gravity waves propagate after the collapse of the initial symmetric peak with respect to axes. Nonlinear interactions create shorter waves propagating around the domain and increasingly more complicated patterns are formed.

\begin{figure*}
\centerline{\includegraphics[width=190pt,height=10pc]{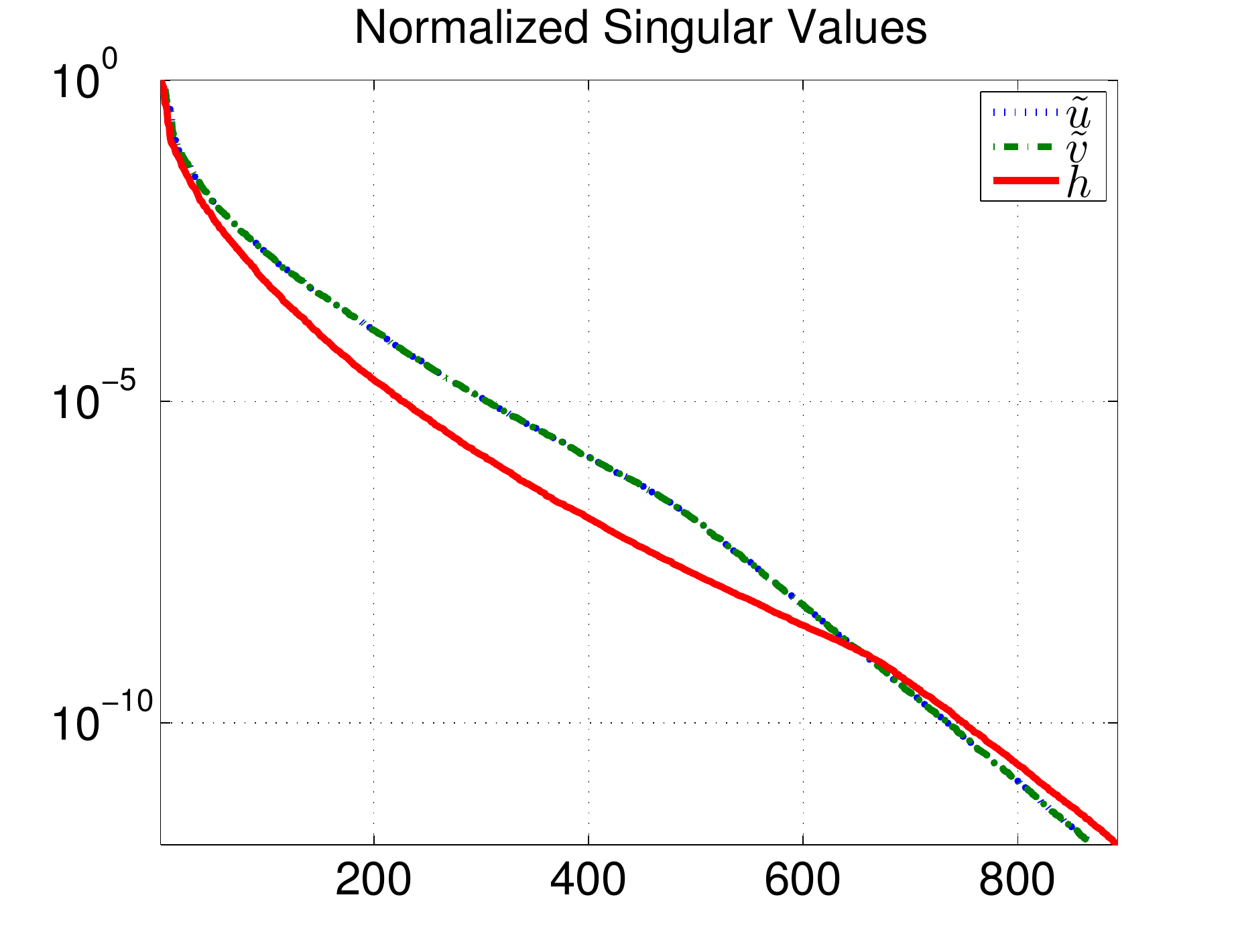}
\includegraphics[width=190pt,height=10pc]{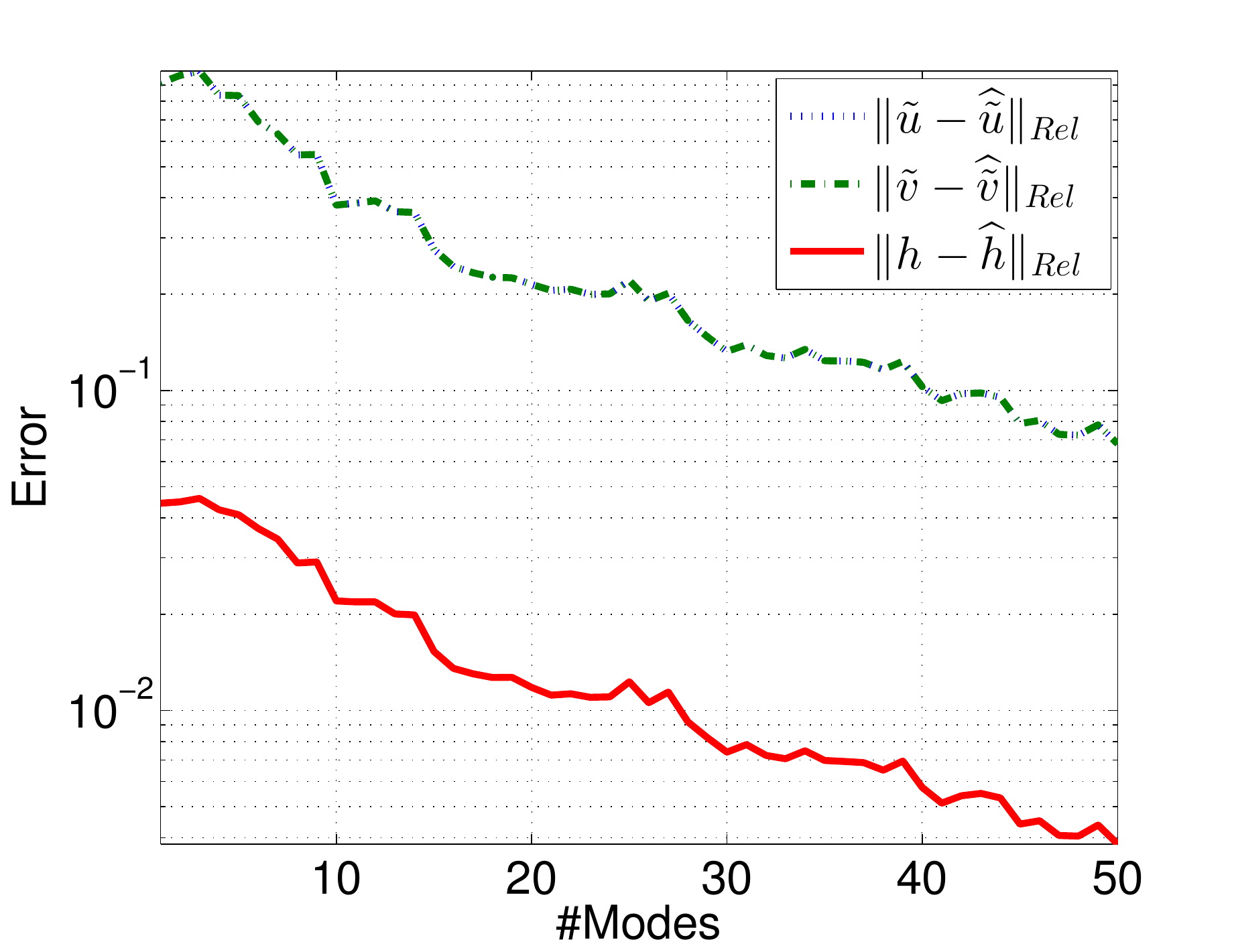}}
\caption{Normalized singular values (left) and relative FOM-ROM errors (right)\label{ex1fig1}}
\end{figure*}

For this test problem, each snapshot matrix $S_{\tilde{u}}$, $S_{\tilde{v}}$ and $S_h$ has sizes $10000\times 1000$.
According to the energy criteria \eqref{energy_criteria}, we take $n=30$ POD modes and $m=200$ DEIM modes.
In Figure~\ref{ex1fig1}, left, the singular values decay slowly for each component, which is the characteristic of the problems with wave phenomena in fluid dynamics \cite{Ohlberger16}. Due to the slow decay of the singular values,  FOM-ROM errors for all components with varying number of POD modes in Figure~\ref{ex1fig1}, right, decrease slowly with small oscillations.

\begin{figure*}
\centerline{\includegraphics[width=190pt,height=10pc]{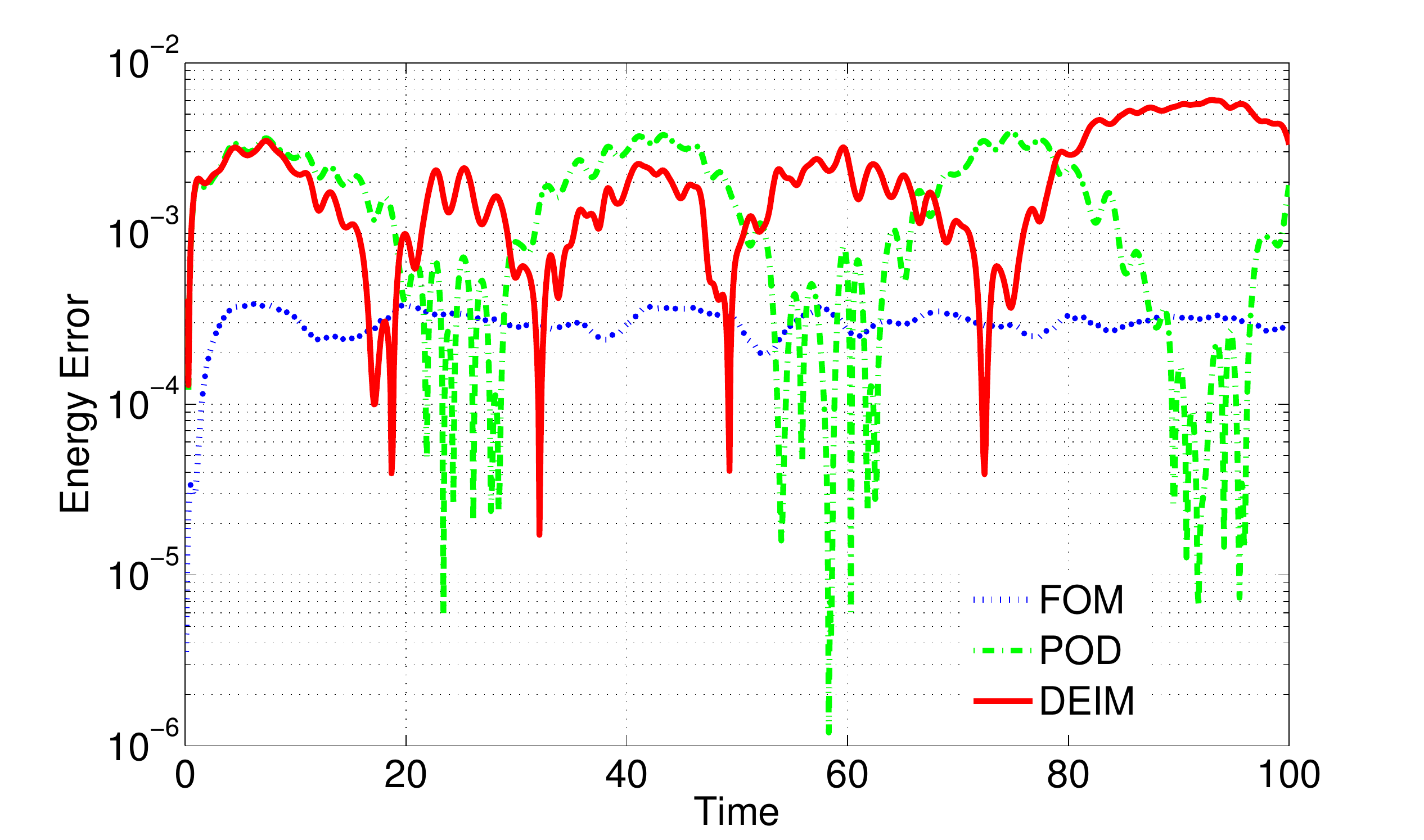}
\includegraphics[width=190pt,height=10pc]{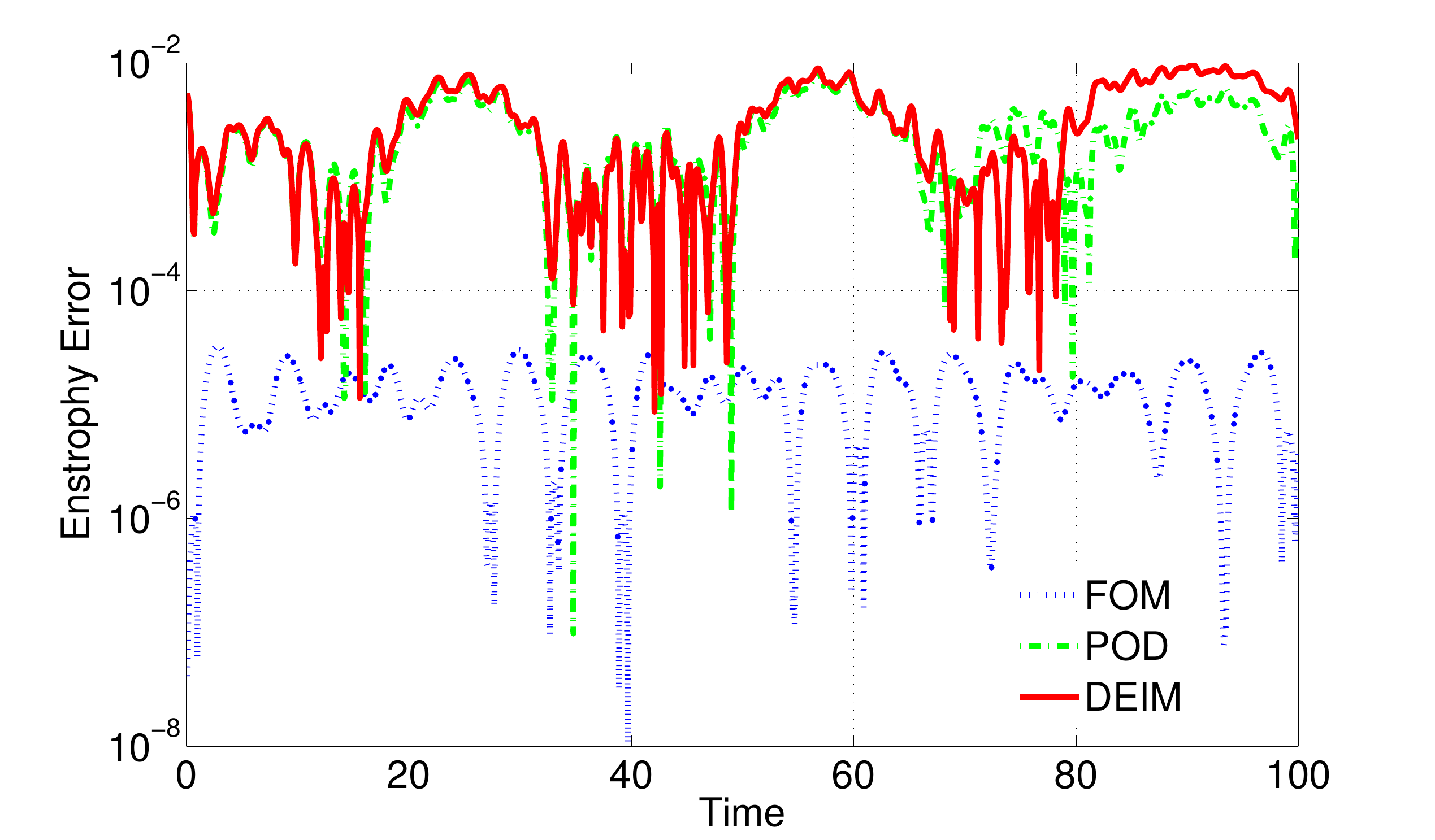}}
\caption{Energy error $|H^k-H^0|$ (left) and enstrophy error $|Z^k-Z^0|$ (right)\label{ex1fig2}}
\end{figure*}

\begin{figure*}
\centerline{\includegraphics[width=130pt,height=9pc]{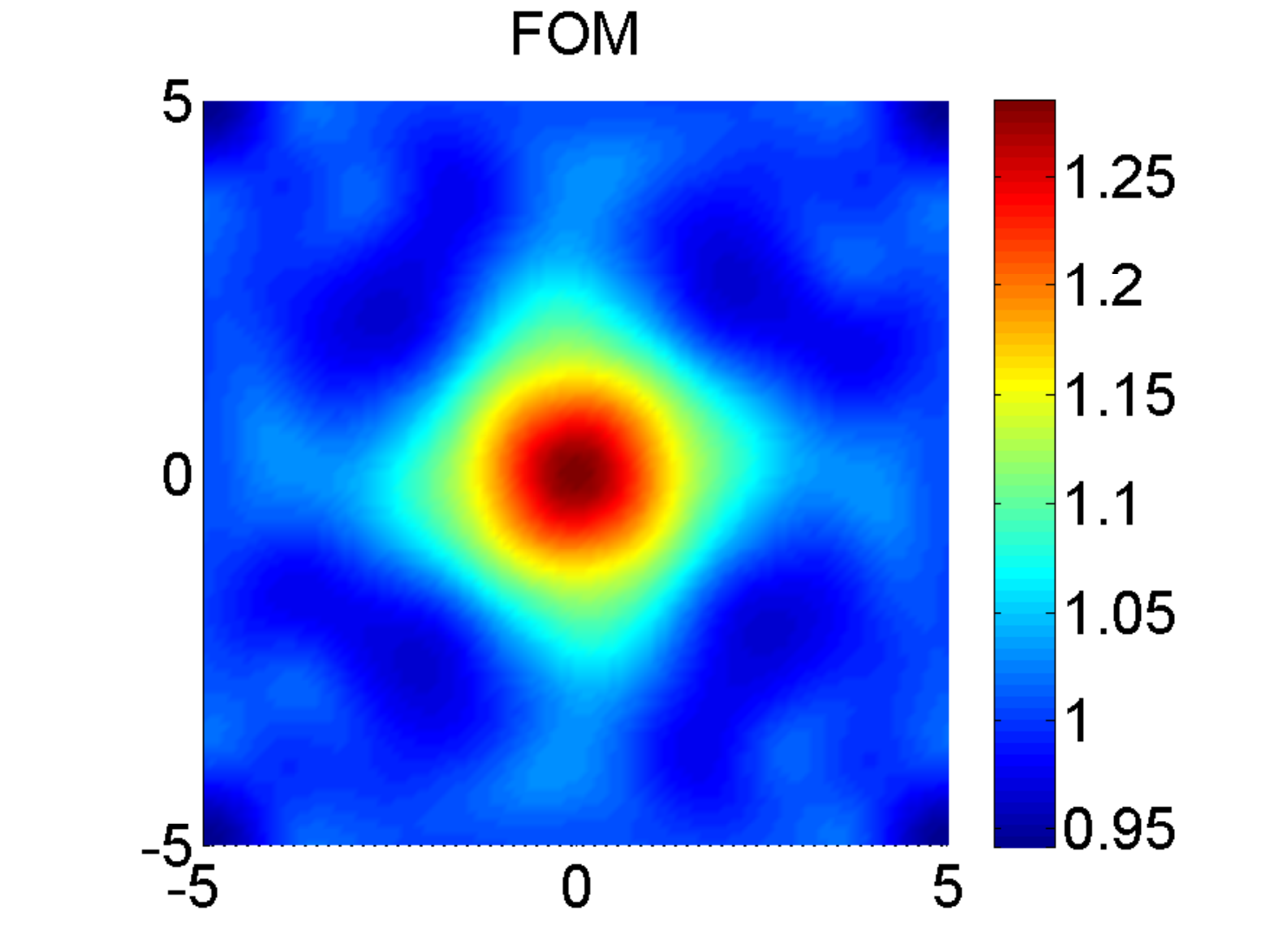}
\includegraphics[width=130pt,height=9pc]{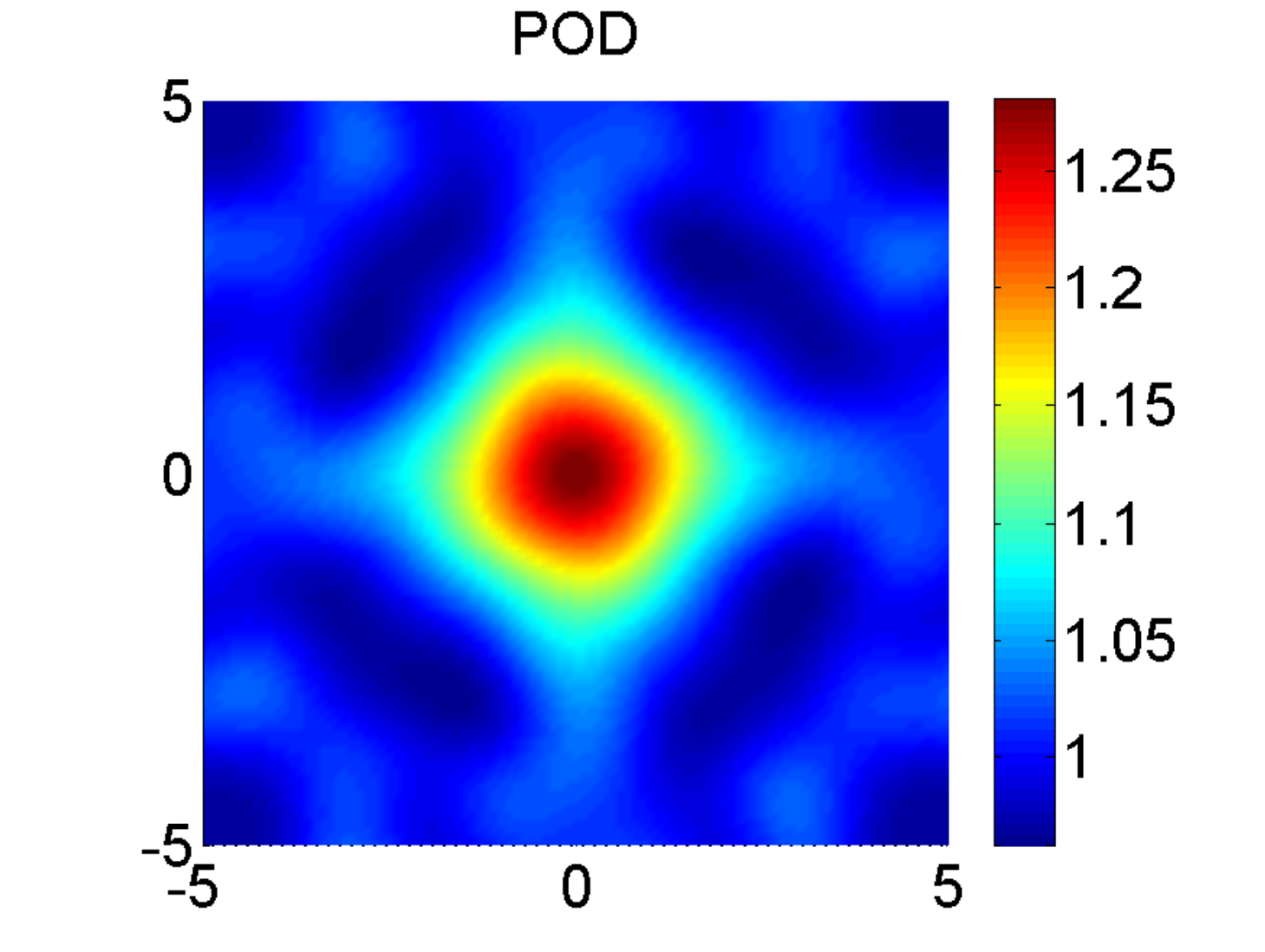}
\includegraphics[width=130pt,height=9pc]{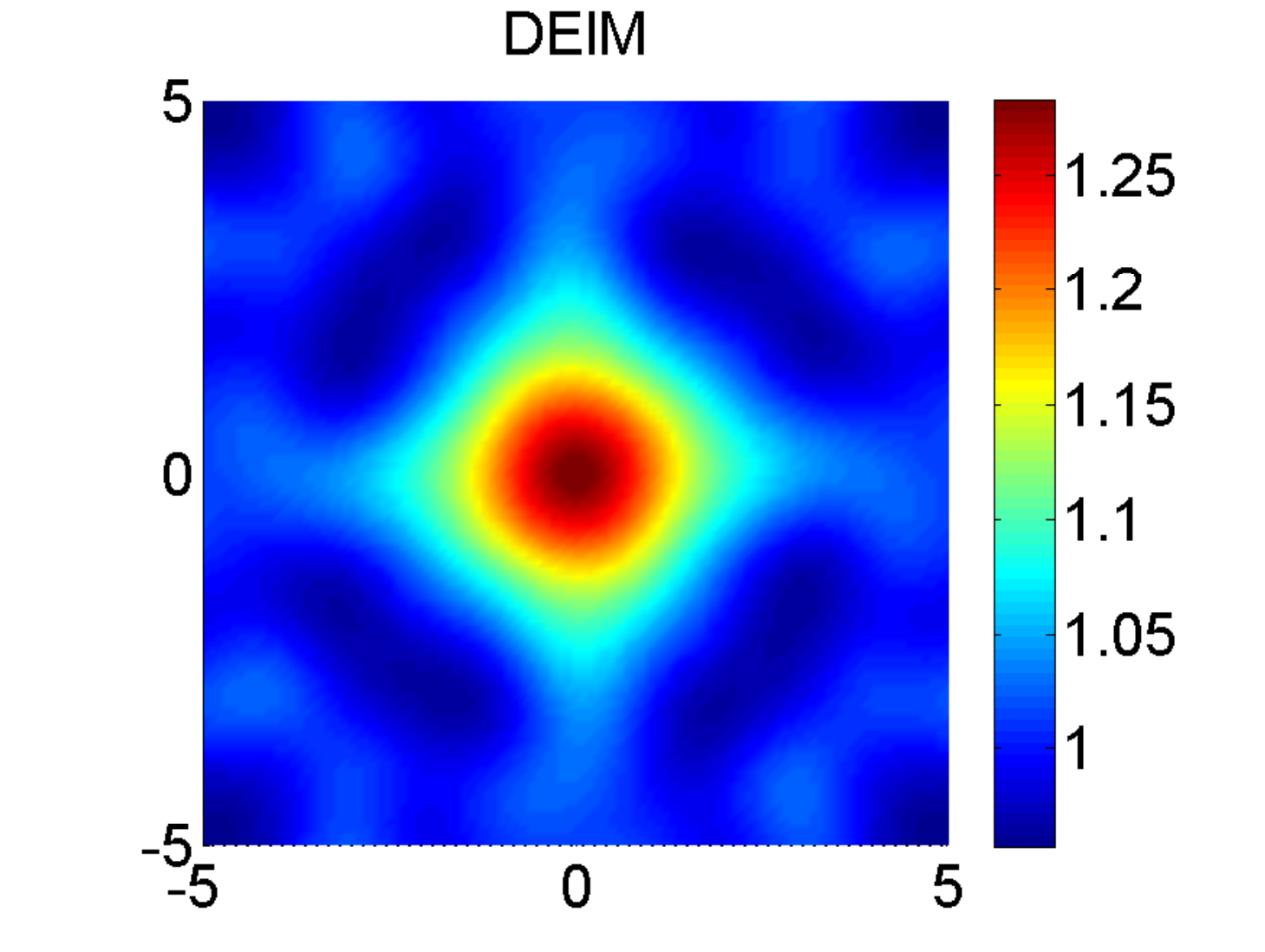}}
\caption{Full and reduced solutions for the height $\bm{h}$ at the final time\label{ex1fig3}}
\end{figure*}

\begin{figure*}
\centerline{\includegraphics[width=130pt,height=9pc]{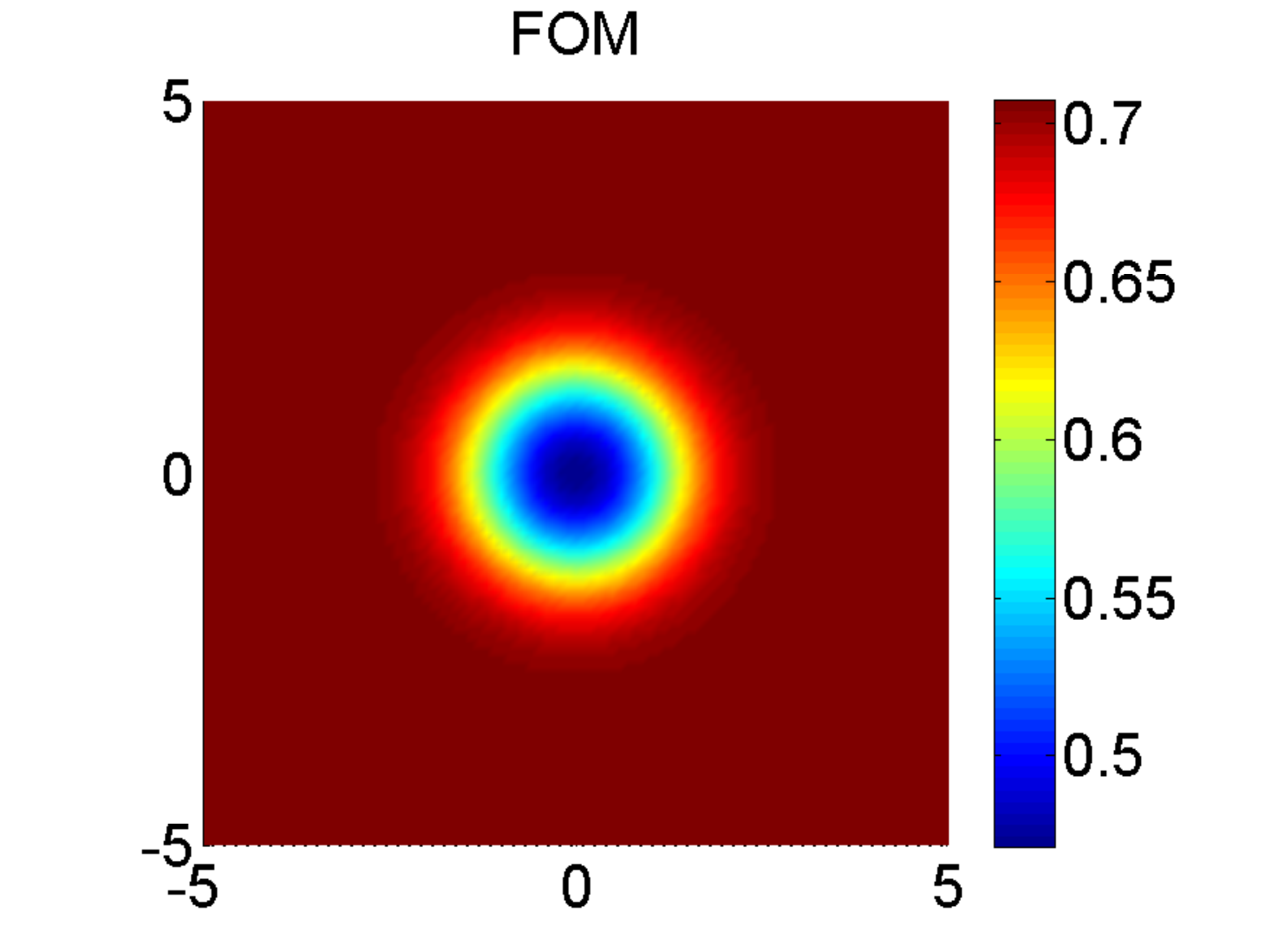}
\includegraphics[width=130pt,height=9pc]{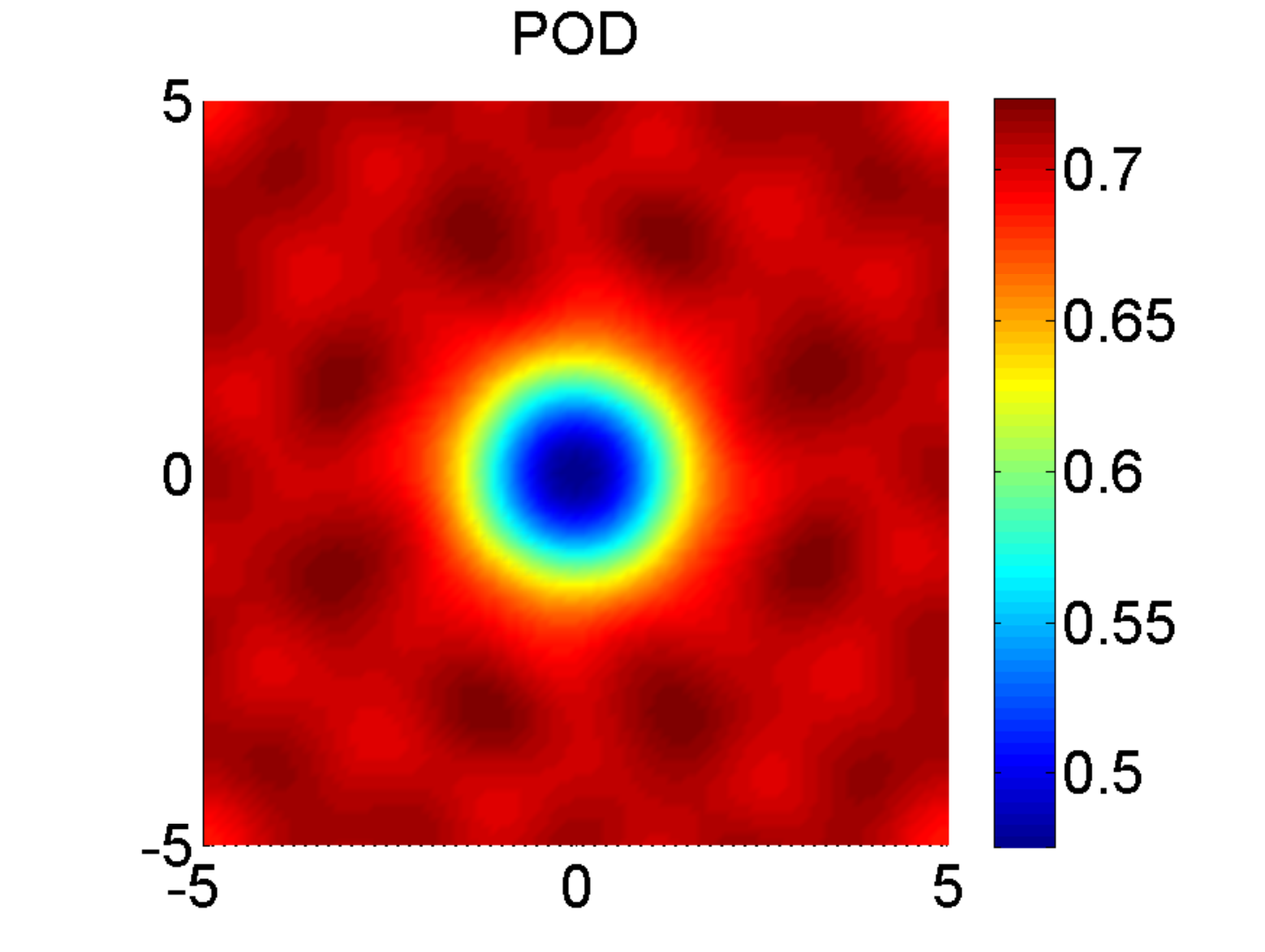}
\includegraphics[width=130pt,height=9pc]{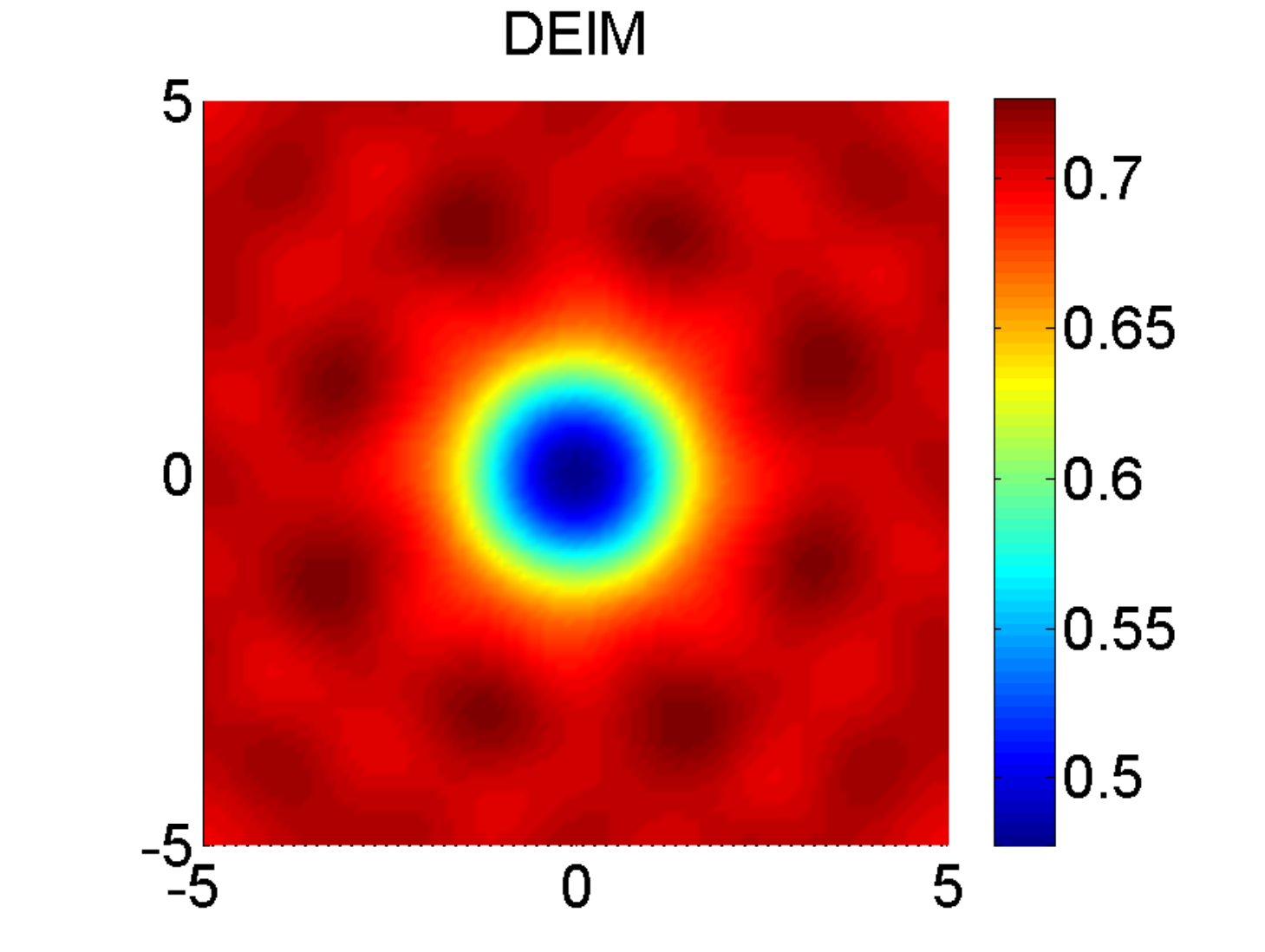}}
\caption{Full and reduced solutions for the potential vorticity $\bm{q}$ at the final time\label{ex1fig5}}
\end{figure*}

The energy and the enstrophy errors in Figure~\ref{ex1fig2} show small drifts with bounded oscillations over the time, i.e. they are preserved approximately at the same level of accuracy.
In Figures~\ref{ex1fig3}-\ref{ex1fig5}, the height $\bm{h}$ and the potential vorticity $\bm{q}$ are shown at the final time. It is seen from the Figures  ~\ref{ex1fig3}-\ref{ex1fig5} and Tables~\ref{tbl1}-\ref{tbl2} that reduced solutions, conserved reduced quantities are of an acceptable level of accuracy. The speed-up factors in Table~\ref{tbl3} shows that the ROM with DEIM increases the computational efficiency further.

\subsection{Single-layer shear instability}
\label{ex2}

We consider the NTSWE on the periodic spatial domain $[0, 10]^2$ and on the time interval $[0,50]$ \cite{Stewart16}.
The initial conditions are given as
\begin{align*}
&h(x, y, 0) = 1+\Delta_h\sin\bigg\{\frac{2\pi}{L}\left[y-\Delta_y \sin\left(\frac{2\pi x}{L}\right)\right]\bigg\},\\
&u(x, y, 0) = -\frac{2\pi\Delta_h}{\Omega^z L}\cos\bigg\{\frac{2\pi}{L}\left[y-\Delta_y \sin\left(\frac{2\pi x}{L}\right)\right]\bigg\},\\
&v(x, y, 0) = -\frac{4\pi^2\Delta_h \Delta_y}{\Omega^z L^2}\cos\bigg\{\frac{2\pi}{L}\left[y-\Delta_y \sin\left(\frac{2\pi x}{L}\right)\right]\bigg\}\cos\left(\frac{2\pi x}{L}\right)
\end{align*}
where $\Delta_h=0.2 $, $ \Delta_y=0.5 $ and the dimensionless spatial domain length $ L=10 $, as the case in the first test example. This problem illustrates the roll-up of an unstable shear layer.

In this test example, each snapshot matrix $S_{\tilde{u}}$, $S_{\tilde{v}}$ and $S_h$ has sizes $10000\times 500$, and the number of POD and DEIM modes are set as $n=18$ and $m=170$, respectively, according to the energy criteria \eqref{energy_criteria}.

\begin{figure*}
\centerline{\includegraphics[width=190pt,height=10pc]{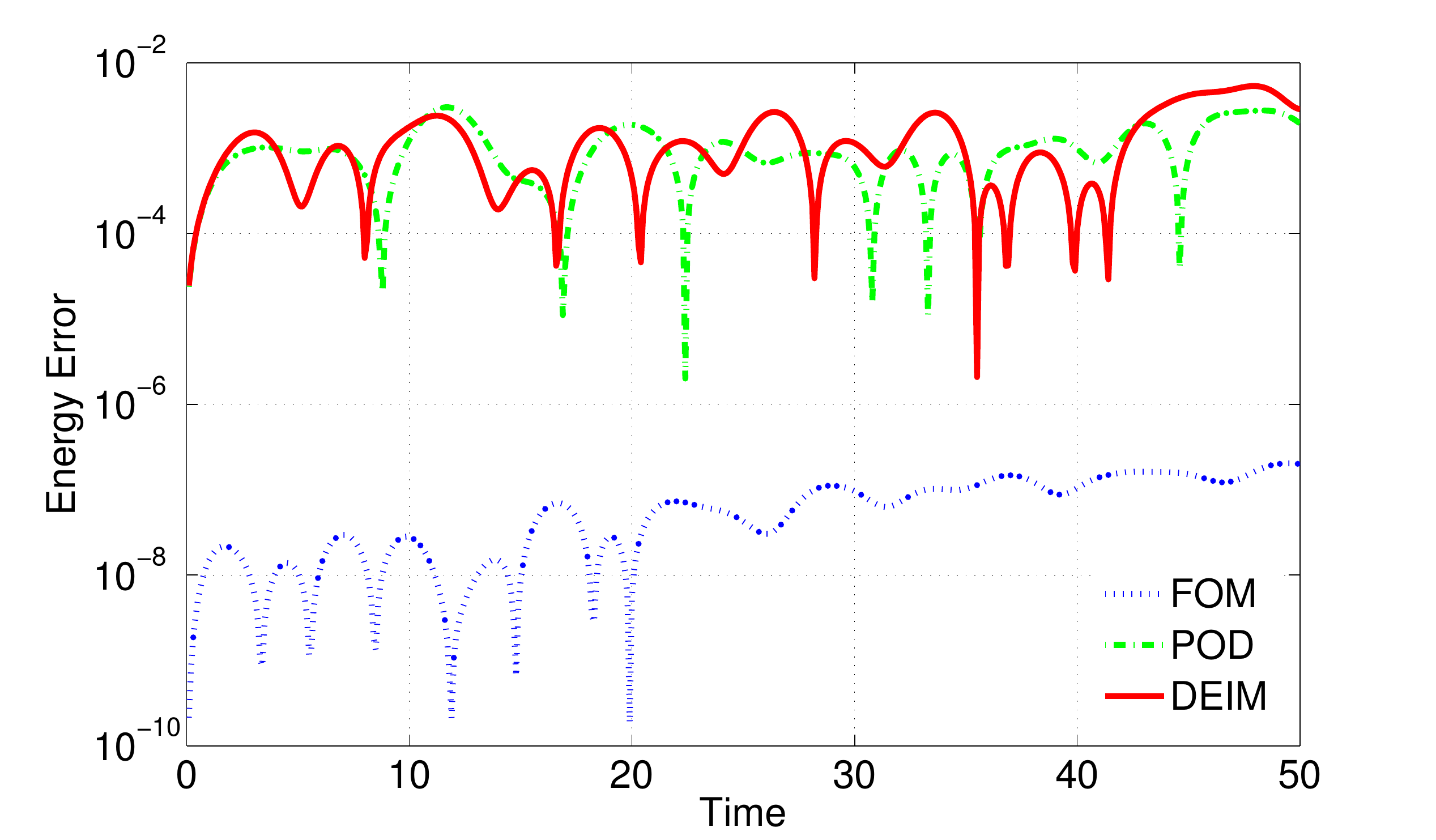}
\includegraphics[width=190pt,height=10pc]{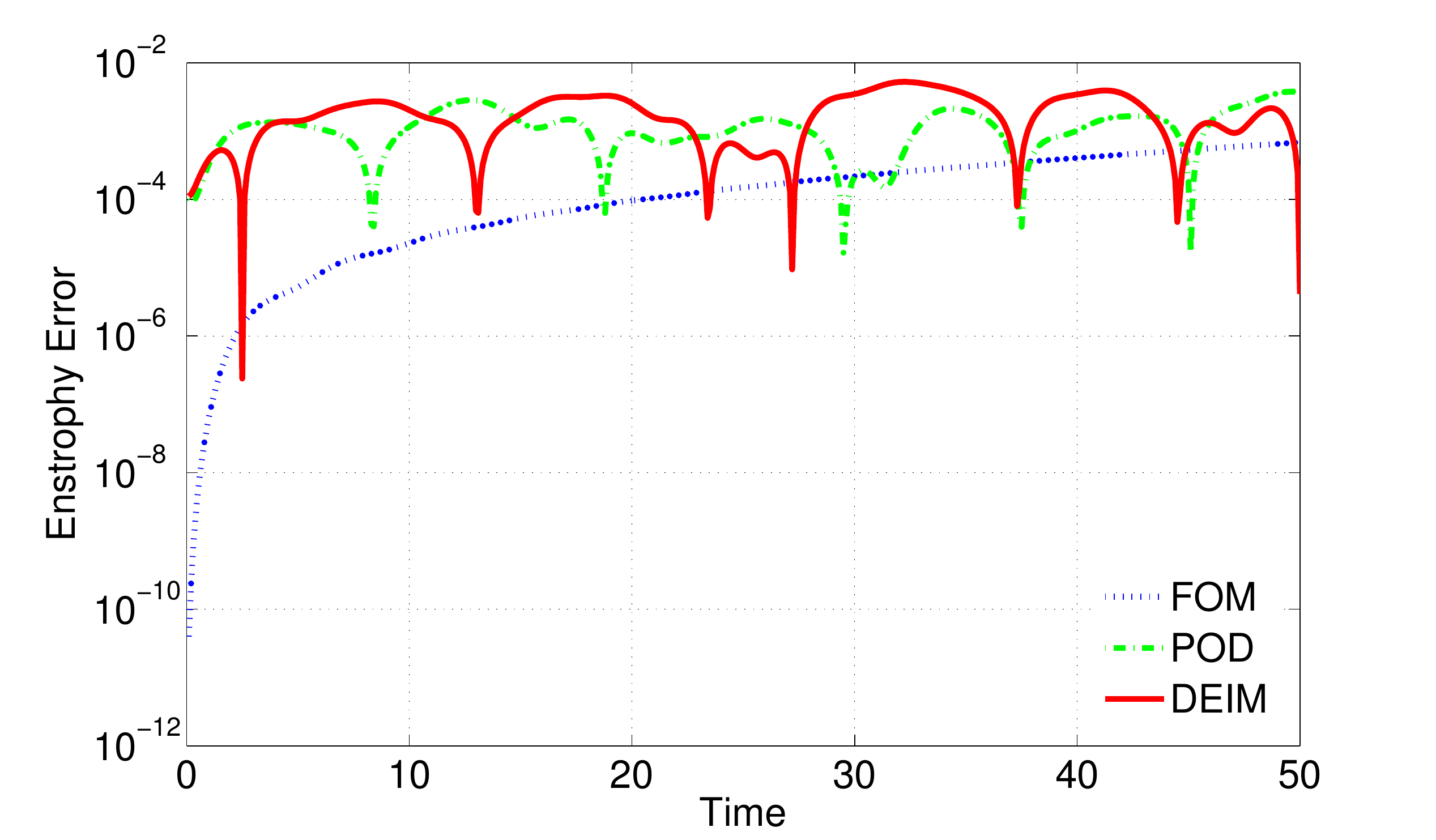}}
\caption{Energy error  $|H^k-H^0|$ (left) and enstrophy error  $|Z^k-Z^0|$ (right)\label{ex2fig2}}
\end{figure*}

\begin{figure*}
\centerline{\includegraphics[width=130pt,height=9pc]{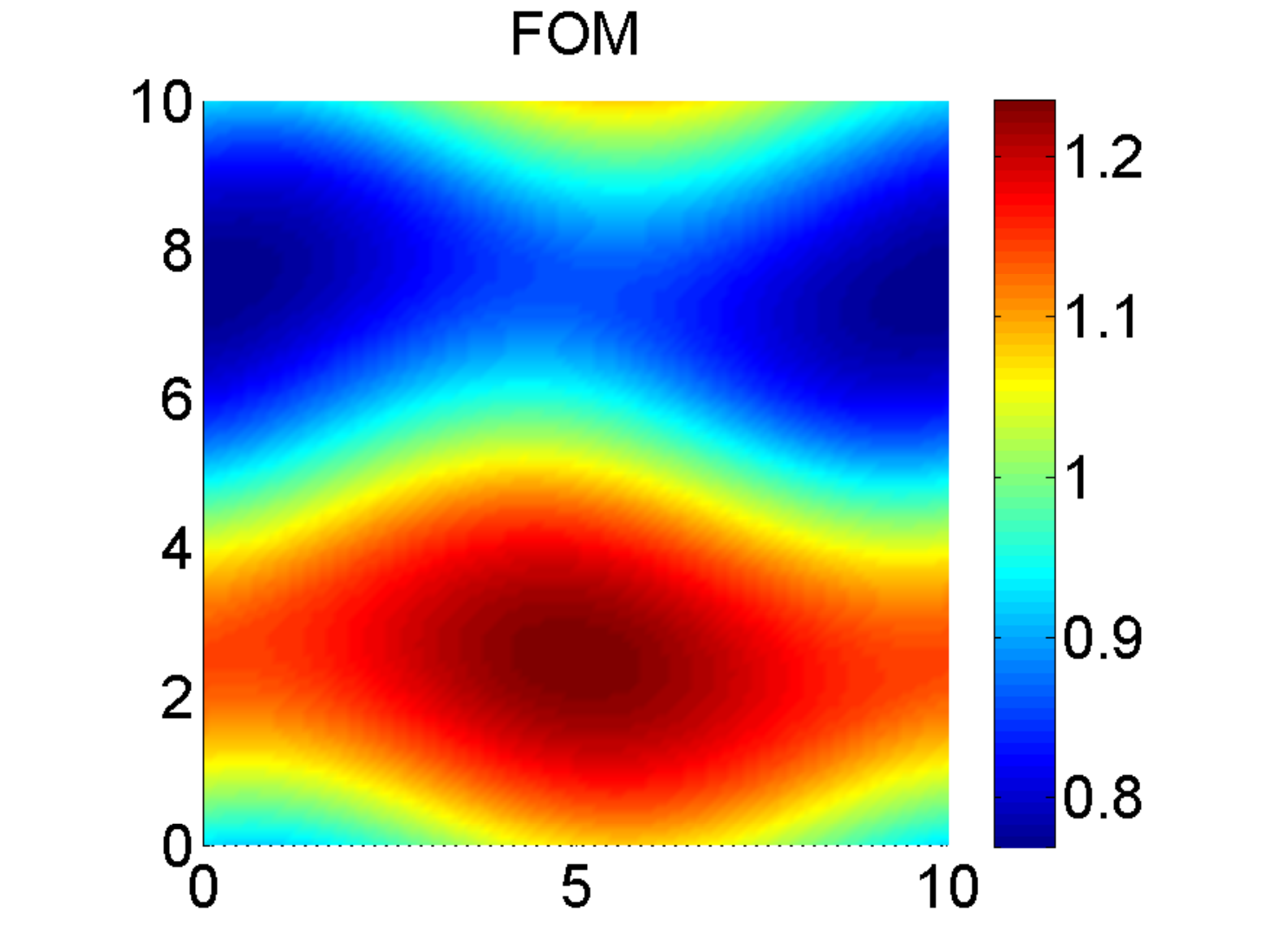}
\includegraphics[width=130pt,height=9pc]{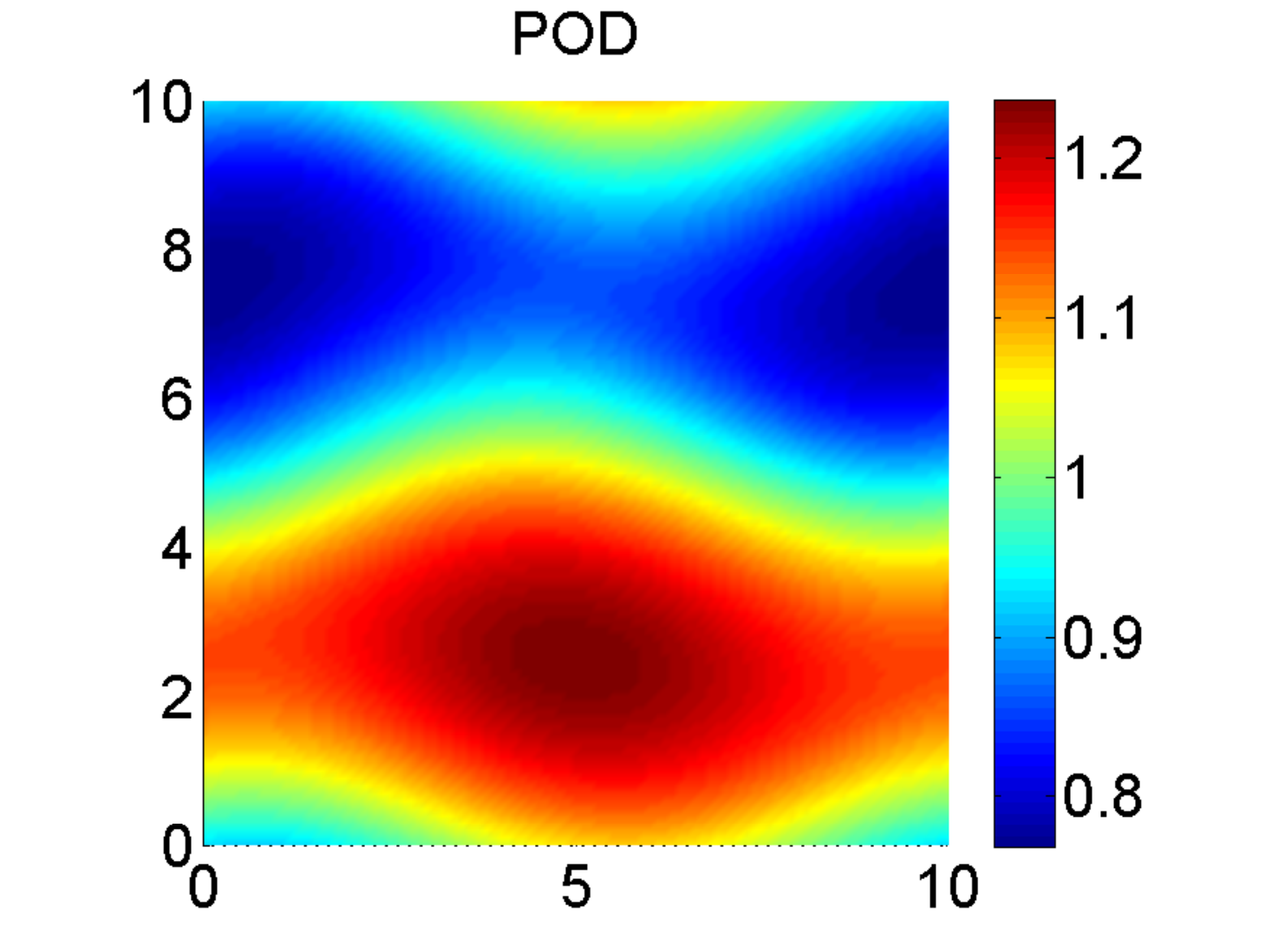}
\includegraphics[width=130pt,height=9pc]{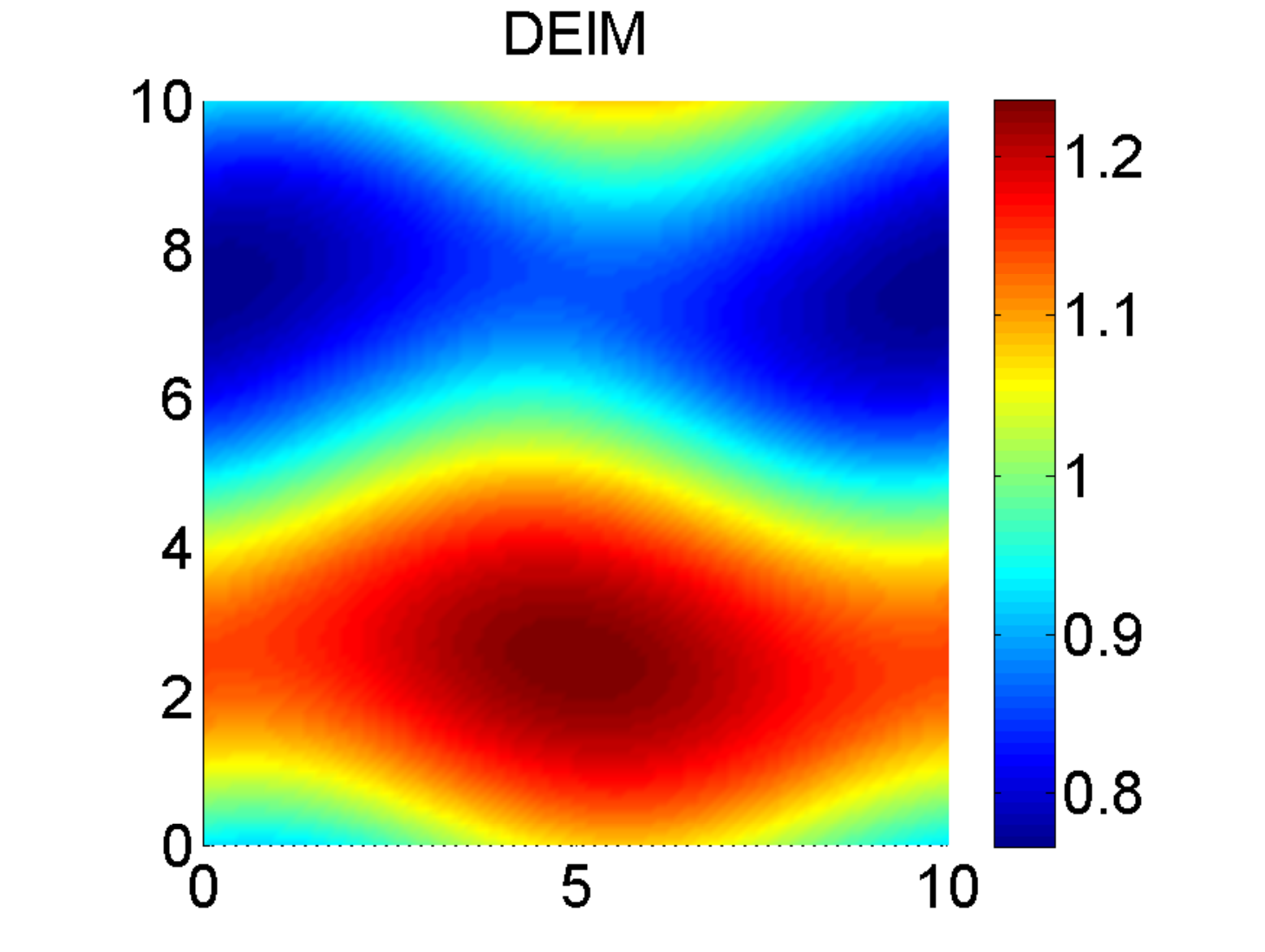}}
\caption{Full and reduced solutions for the height $\bm{h}$ at the final time\label{ex2fig3}}
\end{figure*}

\begin{figure*}
\centerline{\includegraphics[width=130pt,height=9pc]{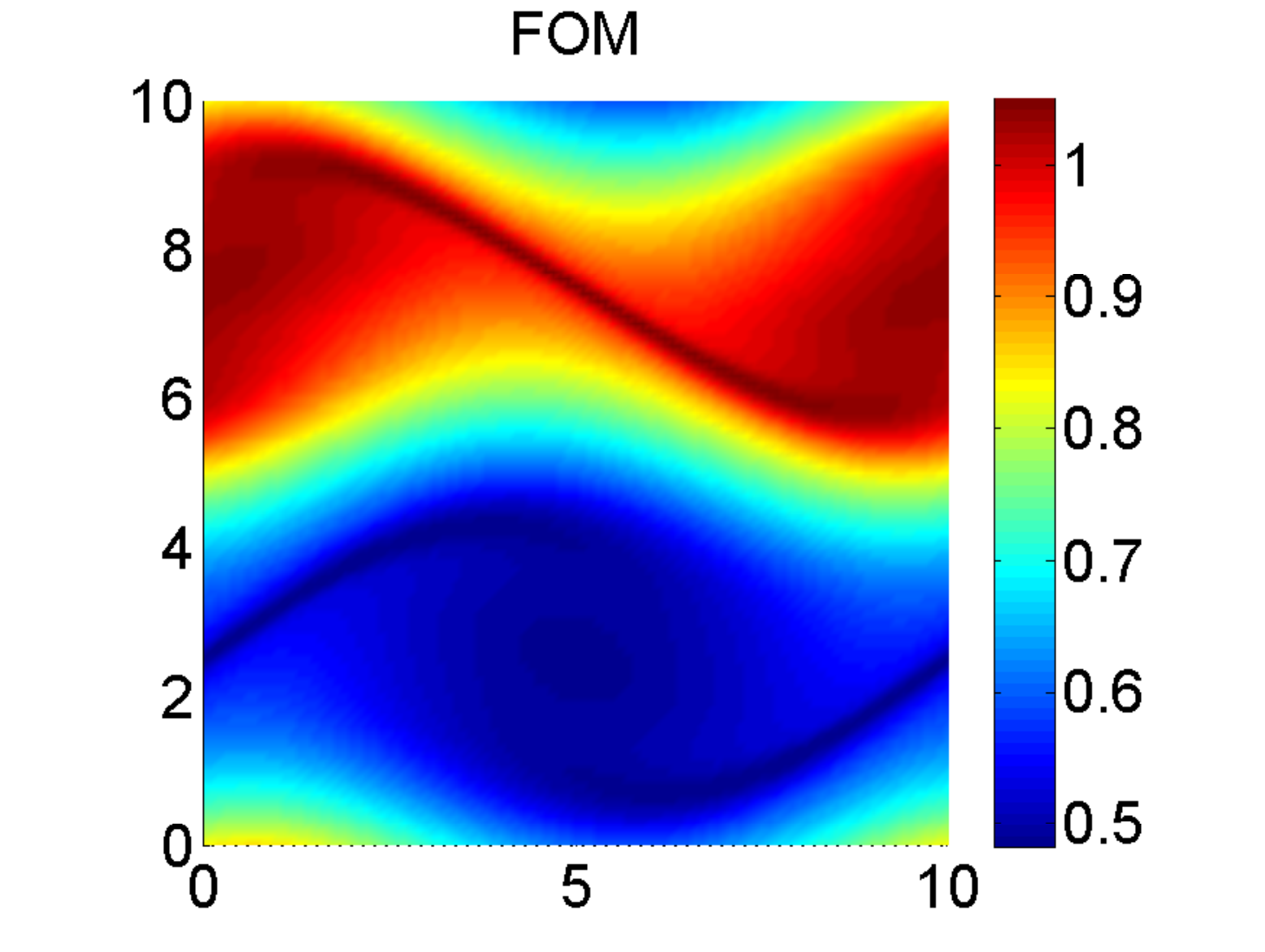}
\includegraphics[width=130pt,height=9pc]{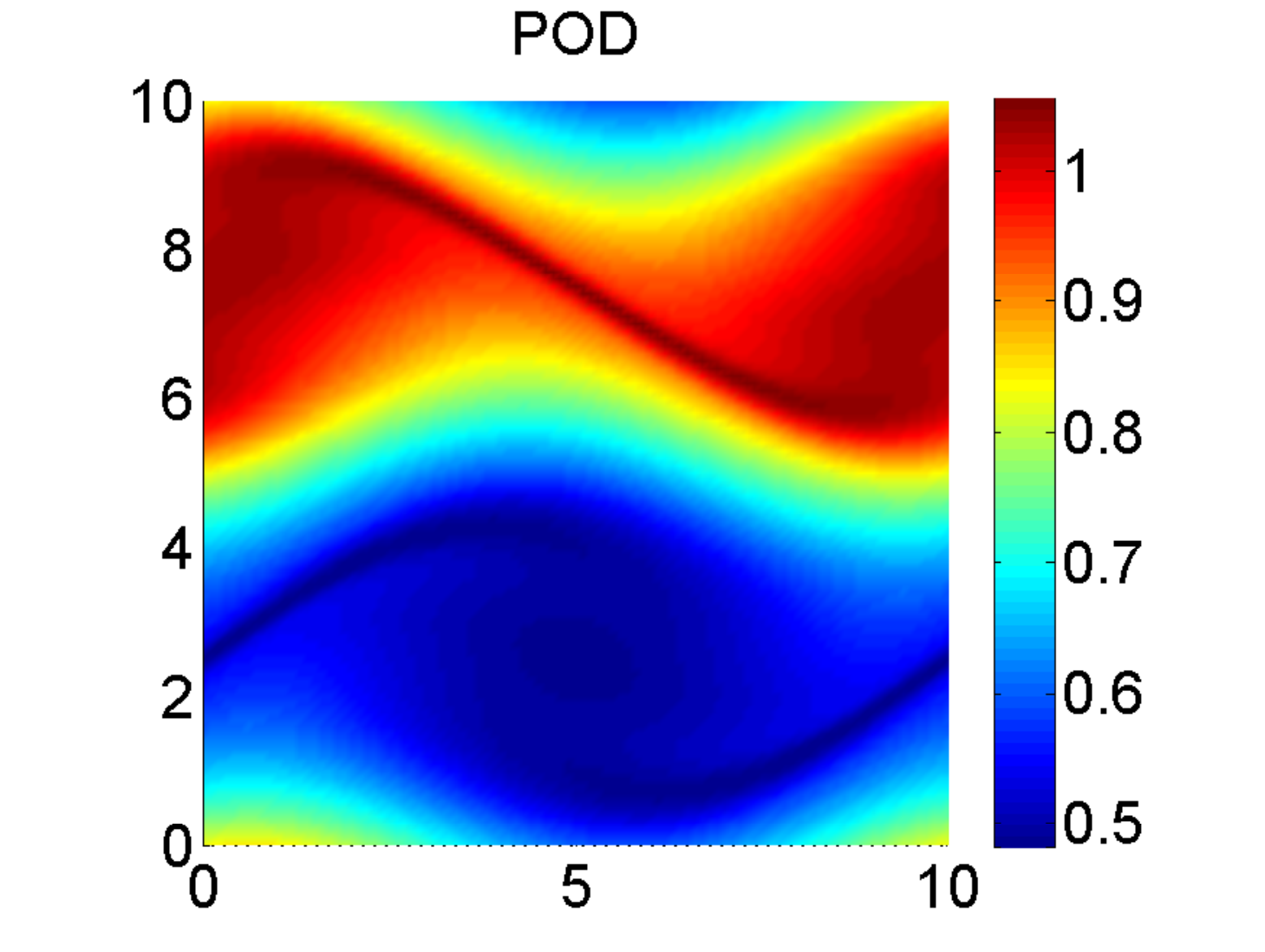}
\includegraphics[width=130pt,height=9pc]{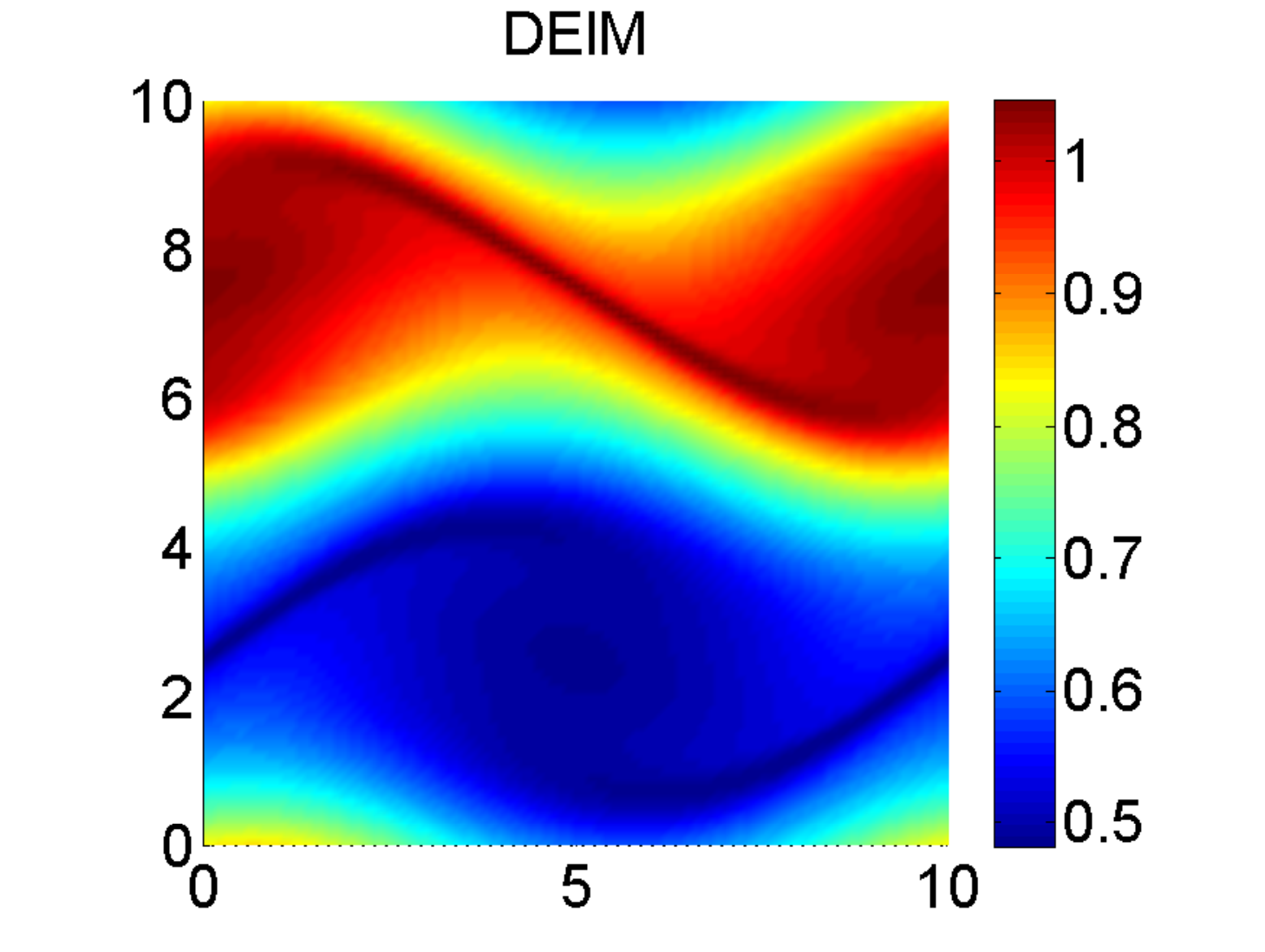}}
\caption{Full and reduced solutions for the potential vorticity $\bm{q}$ at the final time\label{ex2fig5}}
\end{figure*}

The energy and enstropy errors in Figure \ref{ex2fig2} are bounded over time with small oscillations as in the case of the first test example. Similarly, the height $\mathbf{h}$ and  the potential vorticity $\mathbf{q}$  are well approximated by the ROMs  at the final time in Figures~\ref{ex2fig3}-\ref{ex2fig5}. In Tables~\ref{tbl1}-\ref{tbl2} and  Table~\ref{tbl3}, the accuracy and computational efficiency of the reduced approximations are demonstrated.

\begin{table}[htb]
\caption{Time averaged relative $L^2$-errors}%
\label{tbl1}
\begin{tabular}{p{2.5cm}p{2.8cm}p{2.0cm}p{2.0cm}p{2.0cm}}
\hline
 &   & $\|\tilde{\mathbf{u}}-\widehat{\tilde{\mathbf{u}}}\|_{Rel}$  & $\|\tilde{\mathbf{v}}-\widehat{\tilde{\mathbf{v}}}\|_{Rel}$  & $\|\mathbf{h}-\widehat{\mathbf h}\|_{Rel}$ \\
\hline
\multirow{2}{*} {Example~\ref{ex1}} & 30 POD modes  &  1.346e-01 & 1.346e-01  & 7.261e-03   \\
      & 200 DEIM modes  & 1.370e-01  & 1.370e-01  & 7.368e-03   \\
\hline
\multirow{2}{*} {Example~\ref{ex2}} & 18 POD modes  &  2.467e-03 & 9.512e-03  & 2.598e-04   \\
      & 170 DEIM modes  & 3.902e-03  & 1.606e-02  & 4.567e-04   \\
\hline
\end{tabular}
\end{table}

\begin{table}[htb]
\caption{Mean absolute errors between the conserved quantities}%
\label{tbl2}
\begin{tabular}{p{2.5cm}p{3.2cm}p{2.5cm}p{2.5cm}}
\hline
 &   & \textbf{Energy}  & \textbf{Enstrophy} \\
\hline
\multirow{2}{*}{Example~\ref{ex1}} &  30 POD modes  &   1.241e-03 & 1.494e-03    \\
      & 200 DEIM modes & 1.352e-03 & 2.728e-03     \\
\hline
\multirow{2}{*}{Example~\ref{ex2}} &  18 POD modes  &   1.241e-04 & 5.935e-06   \\
      & 170 DEIM modes & 6.458e-04 & 7.119e-04    \\
\hline
\end{tabular}
\end{table}

\begin{table}[htb]%
\caption{CPU time (in seconds) and speed-up factors }
\label{tbl3}
\begin{tabular}{p{1.5cm}p{3cm}p{1.7cm}p{1.7cm}p{1.7cm}p{1.7cm}}
\hline
 &  &  \multicolumn{2}{@{}c@{}}{\textbf{Example~\ref{ex1}}} & \multicolumn{2}{@{}c@{}}{\textbf{Example~\ref{ex2}}} \\
\hline
 &  & CPU time & speed-up & CPU time & speed-up \\
\hline
FOM   &   & 348.3  &   &  244.3 &   \\
\hline
 \multirow{2}{*}{POD}  & basis computation  & 26.2  &  &  8.2 &   \\
	 & online computation  & 233.7  &  1.49  &   101.6 &  2.4 \\
\hline
 \multirow{2}{*}{DEIM}  & basis computation  & 13.8  &   &    3.2 &   \\
	 & online computation  & 23.3  & 14.95 &    16.6 &  14.7 \\
\hline
\end{tabular}
\end{table}

\section{Conclusions}

In contrast to the canonical Hamiltonian systems like the NLS and non-canonical Hamiltonian systems with constant Poisson structure, NTSWE possesses state dependent Poisson structure. In this paper, the  Hamiltonian/energy reduced order modeling approach in \cite{Gong17} is applied by reducing further the computational cost of the ROM in the online stage by exploiting the special structure of the skew-symmetric matrix corresponding to the discretized Poisson structure. The accuracy and computational efficiency of the reduced solutions are demonstrated by numerical examples for the POD and DEIM. Preservation of the energy and enstrophy shows further the stability of the reduced solutions over time. 

\section*{Acknowledgement}
This work was
supported by 100/2000 Ph.D. Scholarship Program of the Turkish Higher Education Council.

%

\end{document}